\documentclass[10pt]{article}
\usepackage{bm,latexsym,amsmath,amsthm,amssymb}
\usepackage{graphicx,indentfirst}
\date{}
\begin{document}
\title{\bf A reaction-diffusion-advection competition model with two free boundaries in heterogeneous time-periodic environment}
\author{ Qiaoling Chen, Fengquan Li\thanks{Corresponding author.\newline
\mbox{}\qquad E-mail: fqli@dlut.edu.cn (F. Li); qiaolingf@126.com (Q. Chen); mathwangfeng@126.com (F. Wang) },
Feng Wang\\
\small School of Mathematical Sciences, Dalian University of Technology,
Dalian 116024, PR China}
\date{}
\maketitle \baselineskip 5pt
\begin{center}
\begin{minipage}{130mm}
{{\bf Abstract.} In this paper, we study the dynamics of a two-species competition model with two different free boundaries in heterogeneous time-periodic environment, where the two species adopt a combination of random movement and advection upward or downward along the resource gradient. We show that the dynamics of this model can be classified into four cases, which forms a spreading-vanishing quartering. The notion of the minimal habitat size for spreading is introduced to determine if species can always spread. Rough estimates of the asymptotic spreading speed of free boundaries and the long time behavior of solutions are also established when spreading occurs. Furthermore, some sufficient conditions for spreading and vanishing are provided.

}

\end{minipage}
\end{center}

\baselineskip=15pt

\section{Introduction}

In this paper, we shall study the dynamical behavior of the solution $(u(x, t), v(x, t), s_1(t), s_2(t))$ to the following reaction-diffusion-advection problem with two free boundaries in the heterogeneous time-periodic environment
\begin{align*}
\left\{\begin{array}{l}
u_t=d_1u_{xx}-\alpha_1 u_x+u(a(x,t)-u-k(x,t)v), \quad 0<x<s_1(t),\quad t>0 ,\\[5pt]
v_t=d_2v_{xx}-\alpha_2 v_x+v(b(x,t)-v-h(x,t)u), \quad 0<x<s_2(t),\quad t>0 ,\\[5pt]
u_{x}(0, t)=v_{x}(0, t)=0, \quad t>0,\\[5pt]
u\equiv 0,\quad x\geq s_1(t),~ t>0;\quad v\equiv 0,\quad x\geq s_2(t),~ t>0,\\[5pt]
s_1'(t)=-\mu_1 u_{x}(s_1(t), t), \quad t>0;\quad s_2'(t)=-\mu_2 v_{x}(s_2(t), t), \quad t>0,\\[5pt]
s_1(0)=s_1^0,~s_2(0)=s_2^0,~ u(x, 0)=u_0(x),~v(x, 0)=v_0(x), \quad 0\leq x< \infty,
\end{array}\right.
 \tag{1.1}
\end{align*}
where $u(x,t)$ and $v(x,t)$ represent the population densities of two competing species at the position $x$ and time $t$; $d_1, d_2$ and $\alpha_1, \alpha_2$ are random diffusion rates and advection rates of species $u, v$, respectively; $\mu_1, \mu_2$ measure the intention to spread into new territories of $u, v$, respectively; $a(x, t), b(x,t)$ and $k(x, t), h(x,t)$ are the intrinsic growth and crowding strength of species $u, v$, respectively, which satisfy the following conditions
\begin{align*}
\left\{\begin{array}{l}
(i)\quad a, b, k, h\in C^{\nu_0, \frac{\nu_0}{2}}([0, \infty)\times\mathbb{R})~
\mbox{for some}~ \nu_0\in(0, 1)~\mbox{and are T-periodic in time t for some}~ T>0;\\[5pt]
(ii)\quad \mbox{there are positive H\"{o}lder continuous and T-periodic functions}~
a_*(t), a^*(t), b_*(t), b^*(t), k_*(t),\\[5pt]
k^*(t), h_*(t)~\mbox{and}~h^*(t)~\mbox{such that}~a_*(t)\leq a(x, t)\leq a^*(t),~b_*(t)\leq b(x, t)\leq b^*(t),~k_*(t)\leq k(x, t)\leq k^*(t),~\\[5pt]
h_*(t)\leq h(x, t)\leq h^*(t),~\mbox{for all}~ (x, t)\in [0, \infty)\times[0,T].
\end{array}\right.
\tag{H1}
\end{align*}
All the parameters are assumed to be positive and the initial data $(u_0, v_0, s_1^0, s_2^0)$ satisfy
\begin{align*}
\left\{\begin{array}{l}
u_0\in C^2([0, s_1^0]),\quad v_0\in C^2([0, s_2^0]),\quad u_0'(0)=v_0'(0)=0,\quad s_1^0>0, ~ s_2^0>0,\\[5pt]
u_0(x)>0~\mbox{for}~ x\in[0, s_1^0),\quad u_0(x)=0~\mbox{for}~ x\geq s_1^0,\\[5pt]
v_0(x)>0~\mbox{for}~ x\in[0, s_2^0),\quad v_0(x)=0~\mbox{for}~ x\geq s_2^0.
\end{array}\right.
 \tag{1.2}
\end{align*}

Ecologically, problem $(1.1)$ may be viewed as a model describing the spreading of two competing species along the same direction with two different free boundaries $x=s_1(t)$ and $x=s_2(t)$,
which may intersect each other at some time, in the heterogeneous time-periodic environment. For simplicity, the left boundary is fixed such that no flux can across the left boundary. We assume that the species can only spread further into the environment from the right end of the initial region, and the spreading front expands
at a speed that is proportional to the population gradient at the front, which gives rise
to the Stefan conditions $s_1'(t)=-\mu_1 u_{x}(s_1(t), t)$ and $s_2'(t)=-\mu_2 v_{x}(s_2(t), t)$.

In the absence of one species $v$, namely $v\equiv0$, problem $(1.1)$ reduces to a reaction-diffusion-advection model with a free boundary in heterogeneous time-periodic environment, which will be considered in Section 3 of this paper later.
For the non-periodic case (i.e. $a$ is independent of time $t$), \cite{gll13,gll14} has studied the dynamics of this problem and \cite{glz15,zyw15,kam14} investigated a more general case, in which  the reaction term $u(a-u)$ is replaced by $f(u)$ including monostable, bistable and combustion type. Similar work but for a SIS reaction-diffusion-advection model can be found in \cite{gklz15}. If the effect of advection is ignored (i.e. $\alpha_1=0$), there are many recent results for time-periodic case \cite{clw14,dgp13,w143} and non-periodic case \cite{bdk12,dl10,dg11,dxl13,dbl13,dmz14,k14,pz13,jlz14,llz14,zx14,w142}. In particular, Du and Lin \cite{dl10} are the first ones to study the spreading-vanishing dichotomy of species in the homogeneous environment of dimension one, namely, the species either
spreads successfully or vanishes eventually. Moreover, the asymptotic spreading speed was
established. Here we call that the species $u$ spreads successfully if $s_{1,\infty}:=\lim_{t\rightarrow \infty}s_{1}(t)=+\infty$ and species $u$ persists in the sense that $\liminf_{t\rightarrow\infty}u(\cdot, t)>0$,
and the species $u$ vanishes eventually if $s_{1,\infty}<+\infty$ and $\lim_{t\rightarrow\infty}\|u(\cdot, t)\|_{C([0, s_1(t)])}=0$.

If $\alpha_1=\alpha_2=0$, $(1.1)$ becomes into a diffusive competition problem with two free boundaries in heterogeneous time-periodic environment, which has not been studied before.
For the non-periodic case, there are many different biological considerations to diffusive competition problem associated with $(1.1)$. In \cite{dl13,dwz14,wz15}, the authers studied a competition problem with a free boundary in which an invasive species exits initially in a ball and invades into the environment, while the resident species distributes in the whole space $\mathbb{R}^N$, that is, $s_{2}(t)=\infty$. In \cite{gw12,wz14}, the two weak competition species are assumed to spread along the same boundary, that is, $s_{1}(t)=s_{2}(t)$. In \cite{gw15,chw15}, the authors considered a two-species model with two different free boundaries both for the weak and strong competition case. For the time-periodic case, \cite{w144,clw15} recently studied the special case $s_{1}(t)= s_{2}(t)$ and $s_{2}(t)=\infty$. Similar works but for two-species Lotka-Volterra type predator-prey problems can be found in \cite{zw14,wz14,w14}. We also refer to much earlier works \cite{myy85,l07} in which the environment is assumed to be a bounded domain.

Motivated by the works \cite{gw15,chw15,dgp13,gll14}, we will study the dynamics of problem $(1.1)$ in more generally natural environment including spatial heterogeneity and daily (or seasonal) changes. Small advection terms are introduced to model the dynamical behavior of two competition species.   We will provide a rather complete description of the spreading-vanishing quartering, sharp threshold for spreading and vanishing, rough estimates of the asymptotic spreading speed of the free boundaries and profile of solutions when spreading happens. On the other hand, it is interesting to understand how the spreading of species depends on the initial habitat and system parameters. Inspired by \cite{chw15}, we introduce the notion of the minimal habitat size for spreading to determine whether their spreading can always succeed.

The rest of this paper is organized as follows.
In Section 2, we study a periodic-parabolic eigenvalue problem and give a comparison
principle for $(1.1)$.
In Section 3, we mainly established the spreading-vanishing dichotomy for problem $(1.1)$ with $v\equiv 0$. Our main results of problem $(1.1)$ are given in Section 4, such as the spreading-vanishing quartering, some sufficient conditions for spreading and vanishing, the long-time behavior of solutions and some rough estimates of the asymptotic spreading speed of free boundaries when spreading happens.

\section{Preliminaries}

In this section, we first consider a linear periodic-parabolic eigenvalue problem and then give a comparison principle for $(1.1)$. These results play an important role in later sections.

\subsection{An eigenvalue problem}
Consider the following eigenvalue problem
\begin{align*}
\left\{\begin{array}{l}
\varphi_t=d\varphi_{xx}-\alpha\varphi_x+\gamma(x,t)\varphi+\lambda\varphi, \quad 0<x<l,~0\leq t\leq T,\\[5pt]
\varphi_{x}(0, t)=\varphi(l, t)=0, \quad 0\leq t\leq T,\\[5pt]
\varphi(x, 0)=\varphi(x, T),\quad 0\leq x \leq l.
\end{array}\right.
 \tag{2.1}
\end{align*}

It is well known \cite{hess91} that $(2.1)$ posses a unique principal eigenvalue $\lambda_1(d, \alpha, \gamma, l, T)$, which corresponds to a positive eigenfunction $\varphi\in C^{2,1}([0,l]\times[0,T])$ provided that $\gamma(x,t)$ satisfies the assumption $(H1)$.\\

\noindent\textbf{Lemma 2.1.}
Let $\gamma(x, t)$ be a function satisfying $(H1)$. Then\\
$(i) \lambda_1(d, \alpha, \gamma,\cdot, T)$ is a strictly decreasing continuous function in $[0, +\infty)$ for fixed $d, \alpha, \gamma, T$ and $\lambda_1(d, \alpha, \cdot, l, T)$ is a strictly decreasing
continuous function in the sense that, $\lambda_1(d, \alpha, \gamma_1, l, T)<\lambda_1(d, \alpha, \gamma_2, l, T)$ if the two T-periodic continuous functions $\gamma_1(x,t)$ and $\gamma_2(x, t)$ satisfy $\gamma_1(x, t)\geq \not\equiv \gamma_2(x,t)$ on $[0,l]\times[0, T]$;\\
$(ii) \lambda_1(d, \alpha, \gamma,l,T)\rightarrow+\infty$ as $l\rightarrow 0$;\\
$(iii) \lim_{l\rightarrow+\infty}\lambda_1(d, \alpha, \gamma,l,T)<0$
if $\alpha^{2}<4d\min_{[0, T]}\gamma_{*}(t)$.\\

\noindent\textbf{Proof.} The proof of $(i)$ is similar to that of Lemma 3.2 in \cite{dgp13}. Now we prove $(ii)$ and $(iii)$. It follows from $(i)$ that
\begin{align*}
\lambda_{1}(d,\alpha, \max_{[0,T]}\gamma^{*},l,T)
\leq \lambda_1(d, \alpha, \gamma,l,T)
\leq \lambda_{1}(d,\alpha, \min_{[0,T]}\gamma_{*},l,T).
\end{align*}
It is also easy to see that $\lambda_{1}(d,\alpha, \min_{[0,T]}\gamma_{*},l,T)$ and
$\lambda_{1}(d,\alpha, \max_{[0,T]}\gamma^{*},l,T)$ are the principal eigenvalues of the elliptic problems
\begin{align*}
\left\{\begin{array}{l}
-d\psi_{xx}+\alpha\psi_x-(\min_{[0,T]}\gamma_{*})\psi=\lambda\psi, \quad 0<x<l,\\[5pt]
\psi_{x}(0)=\psi(l)=0\\[5pt]
\end{array}\right.
\end{align*}
and
\begin{align*}
\left\{\begin{array}{l}
-d\psi_{xx}+\alpha\psi_x-(\max_{[0,T]}\gamma^{*})\psi=\lambda\psi, \quad 0<x<l,\\[5pt]
\psi_{x}(0)=\psi(l)=0\\[5pt]
\end{array}\right.
\end{align*}
respectively. It is well known that
\begin{align*}
\lim_{l\rightarrow 0^{+}}\lambda_{1}(d,\alpha, \min_{[0,T]}\gamma_{*},l,T)=+\infty,\quad
\lim_{l\rightarrow +\infty}\lambda_{1}(d,\alpha, \min_{[0,T]}\gamma_{*},l,T)
=\frac{\alpha^{2}}{4d}-\min_{[0,T]}\gamma_{*}
\end{align*}
and
\begin{align*}
\lim_{l\rightarrow 0^{+}}\lambda_{1}(d,\alpha, \max_{[0,T]}\gamma^{*},l,T)=+\infty,\quad
\lim_{l\rightarrow +\infty}\lambda_{1}(d,\alpha, \max_{[0,T]}\gamma^{*},l,T)
=\frac{\alpha^{2}}{4d}-\max_{[0,T]}\gamma^{*},
\end{align*}
which implies $(ii)$ and $(iii)$. \hfill $\Box$

\subsection{A comparison principle}

In this subsection, we give a comparison principle for $(1.1)$. The proof is similar to that of Lemma 2.3 in \cite{chw15}.\\

\noindent\textbf{Lemma 2.2.}
Assume that $T_0\in(0, \infty), \underline{s}_i, \bar{s}_i\in C^1([0,T_0])~(i=1,2), \underline{u}\in C(\overline{D_{T_0}^*})\cap C^{2,1}(D_{T_0}^*)$ with $D_{T_0}^*=\{(x,t):0\leq x\leq \underline{s}_1(t), t\in[0, T_0]\}, \bar{u}\in C(\overline{D_{T_0}^{**}})\cap C^{2,1}(D_{T_0}^{**})$ with $D_{T_0}^{**}=\{(x,t):0\leq x\leq \bar{s}_1(t), t\in[0, T_0]\},
\underline{v}\in C(\overline{E_{T_0}^*})\cap C^{2,1}(E_{T_0}^*)$ with $E_{T_0}^*=\{(x,t):0\leq x\leq \underline{s}_2(t), t\in[0, T_0]\}, \bar{v}\in C(\overline{E_{T_0}^{**}})\cap C^{2,1}(E_{T_0}^{**})$ with $E_{T_0}^{**}=\{(x,t):0\leq x\leq \bar{s}_2(t), t\in[0, T_0]\}$,
and
\begin{align*}
\left\{\begin{array}{l}
\bar{u}_t\geq d_1\bar{u}_{xx}-\alpha_1 \bar{u}_x+\bar{u}(a(x,t)-\bar{u}-k(x,t)\underline{v}), \quad 0<x<\bar{s}_1(t),~0<t\leq T_0,\\[5pt]
\underline{u}_t\leq d_1\underline{u}_{xx}-\alpha_1 \underline{u}_x+\underline{u}(a(x,t)-\underline{u}-k(x,t)\bar{v}),\quad 0<x<\underline{s}_1(t),~0<t\leq T_0,\\[5pt]
\bar{v}_t\geq d_2\bar{v}_{xx}-\alpha_2\bar{v}_x+\bar{v}(b(x,t)-\bar{v}-h(x,t)\underline{u}),\quad 0<x<\bar{s}_2(t), ~0<t\leq T_0,\\[5pt]
\underline{v}_t\leq d_2\underline{v}_{xx}-\alpha_2\underline{v}_x+\underline{v}(b(x,t)-\underline{v}-h(x,t)\bar{u}),\quad 0<x<\underline{s}_2(t),~0<t\leq T_0,\\[5pt]
\bar{u}_x(0,t)=\underline{v}_x(0,t)= 0,~\underline{u}_x(0,t)=\bar{v}_x(0,t)= 0, \quad
0<t\leq T_0,\\[5pt]
\bar{u}\equiv 0,~ x\geq \bar{s}_1(t),~ 0<t\leq T_0;~ \underline{v}\equiv 0,~ x\geq \underline{s}_2(t),~ 0<t\leq T_0,\\[5pt]
\bar{v}\equiv 0,~ x\geq \bar{s}_2(t),~ 0<t\leq T_0;~ \underline{u}\equiv 0,~ x\geq \underline{s}_1(t),~ 0<t\leq T_0,\\[5pt]
\bar{s}_1^{\prime}(t)\geq-\mu_1\bar{u}_x(\bar{s}_1(t),t), ~ \underline{s}_2^{\prime}(t)\leq-\mu_2\underline{v}_x(\underline{s}_2(t),t),\quad 0<t\leq T_0,\\[5pt]
\bar{s}_1(0)\geq s_{1}^0,~ \underline{s}_2(0)\leq s_{2}^0,\\[5pt]
\bar{u}(x,0)\geq u_0(x),\quad 0\leq x\leq s_{1}^0,\\[5pt]
\underline{v}(x,0)\leq v_0(x),\quad 0\leq x\leq \underline{s}_{2}^0.
\end{array}\right.
\end{align*}
Then the solution $(u,v,s_1,s_2)$ of $(1.1)$ satisfies
\begin{align*}
\bar{s}_1(t)\geq s_1(t), \quad \underline{s}_2(t)\leq s_2(t),\quad u(x,t)\leq \bar{u}(x,t), \quad v(x,t)\geq \underline{v}(x,t),~\mbox{for}~(x,t)\in[0, \infty)\times[0,T_0],\\[5pt]
\bar{s}_2(t)\geq s_2(t), \quad \underline{s}_1(t)\leq s_1(t),\quad u(x,t)\geq \underline{u}(x,t), \quad \bar{v}(x,t)\geq v(x,t),~\mbox{for}~(x,t)\in[0, \infty)\times[0,T_0].
\end{align*}

In the following, we shall use $u^{\mu_1, \mu_2}, v^{\mu_1, \mu_2}, s_1^{\mu_1, \mu_2}$ and $s_2^{\mu_1, \mu_2}$ to emphasize the dependence of the solutions on $\mu_1$ and $\mu_2$.
By Lemma 2.2, we have the following result.\\

\noindent\textbf{Corollary 2.1.}
Let $(u_{0}, v_{0}, s_{1}^{0}, s_{2}^{0})$ and all parameters in $(1.1)$ be fixed except for
the parameter $\mu_{1}$ and $\mu_{2}$. If $\mu_{1}\leq \bar{\mu}_{1}$ and $\mu_{2}\geq \underline{\mu}_{2}$, then
\begin{align*}
&s_{1}^{\bar{\mu}_{1},\underline{\mu}_{2}}(t)\geq s_{1}^{\mu_{1},\mu_{2}}(t)
\quad \mbox{and}\quad s_{2}^{\bar{\mu}_{1},\underline{\mu}_{2}}(t)\leq s_{2}^{\mu_{1},\mu_{2}}(t)
\quad \mbox{for}~t>0;\\
&u^{\bar{\mu}_{1},\underline{\mu}_{2}}(x,t)\geq u^{\mu_{1},\mu_{2}}(x,t)
\quad \mbox{for}~x\in(0, s_{1}^{\mu_{1},\mu_{2}}(t)),~t>0;\\
&v^{\bar{\mu}_{1},\underline{\mu}_{2}}(x,t)\leq v^{\mu_{1},\mu_{2}}(x,t)
\quad \mbox{for}~x\in(0, s_{2}^{\bar{\mu}_{1},\underline{\mu}_{2}}(t)),~t>0.
\end{align*}

\section{The spreading-vanishing dichotomy for $(1.1)$ with $v\equiv 0$}

In this section, we mainly established the spreading-vanishing dichotomy for problem
$(1.1)$ with $v\equiv 0$, that is, a reaction-diffusion-advection model with a free boundary
in time-periodic environment.
To achieve it, we need to consider the existence and uniqueness of positive bounded solution of a
$T$-periodic boundary-value problem in the half line.

\subsection{Positive bounded solutions of a $T$-periodic boundary-value problem in the half line}

In this subsection, we first consider the following initial-boundary value problem in the half line
\begin{align*}
\left\{\begin{array}{l}
p_t=dp_{xx}-\alpha p_x+p(\gamma(x,t)-p), \quad 0<x<\infty,~t>0,\\[5pt]
p_{x}(0, t)=0, \quad t>0,\\[5pt]
p(x, 0)=p_0(x),\quad 0\leq x < \infty,
\end{array}\right.
 \tag{3.1}
\end{align*}
and the corresponding T-periodic boundary value problem
\begin{align*}
\left\{\begin{array}{l}
P_t=dP_{xx}-\alpha P_x+P(\gamma(x,t)-P), \quad 0<x<\infty,~0\leq t\leq T,\\[5pt]
P_{x}(0, t)=0, \quad 0\leq t\leq T,\\[5pt]
P(x, 0)=P(x, T),\quad 0\leq x <\infty,
\end{array}\right.
 \tag{3.2}
\end{align*}
where $\gamma(x,t)$ is a $T$-periodic function satisfying the assumption $(H1)$ and $p_0(x)$ is a bounded nontrivial and nonnegative continuous function.

By $(H1)$, we have that
\begin{align*}
\gamma^{\infty}(t):=\limsup_{x\rightarrow+\infty}\gamma(x,t)\leq \gamma^*(t),\quad
\gamma_{\infty}(t):=\liminf_{x\rightarrow+\infty}\gamma(x,t)\geq \gamma_*(t),
\end{align*}
and $\gamma_{\infty}$ and $\gamma^{\infty}$ are T-periodic functions. We assume that these functions are H\"{o}lder continuous.\\

\noindent\textbf{Proposition 3.1.}
Let $\gamma(x, t)$ be a function satisfying $(H1)$ and
$0<\alpha<2\sqrt{d\min_{[0,T]}\gamma_{*}(t)}$, then $(3.2)$
has a unique positive bounded solution $P(x,t)\in C^{2,1}([0,\infty)\times[0, T])$. Furthermore,
\begin{align*}
\mathfrak{p}_{\infty}(t)\leq \liminf_{x\rightarrow \infty}P(x,t)
\leq\limsup_{x\rightarrow \infty}P(x,t)\leq \mathfrak{p}^{\infty}(t)
\end{align*}
uniformly in $[0, T]$, where $\mathfrak{p}_{\infty}(t), \mathfrak{p}^{\infty}(t)$ are the unique positive solutions of the following $T$-periodic ordinary differential problems:
\begin{align*}
\mathfrak{p}^{\prime}(t)=\mathfrak{p}(\gamma_{\infty}(t)-\mathfrak{p}),
\quad \mathfrak{p}(0)=\mathfrak{p}(T),
\tag{3.3}
\end{align*}
and
\begin{align*}
\mathfrak{p}^{\prime}(t)=\mathfrak{p}(\gamma^{\infty}(t)-\mathfrak{p}),
\quad \mathfrak{p}(0)=\mathfrak{p}(T),
\tag{3.4}
\end{align*}
respectively. Moreover, let $p(x,t)$ be the positive solution of $(3.1)$, then $\lim_{n\rightarrow\infty}p(x,t+nT)=P(x,t)$ locally uniformly in $[0, \infty)$.\\

\noindent\textbf{Proof.}
In \cite{pw12}, the authors obtained a similar result for a diffusive logistic equation in $\mathbb{R}^{N}$. Since $(3.2)$ contains a advection term and is considered in the half line, we give the details of proof. We divide our proof in several steps.

\textbf{Step 1.} the existence of the minimal positive solution of $(3.2)$.

Since $\gamma$ satisfies $(H1)$ and
$0<\alpha<2\sqrt{d\min_{[0,T]}\gamma_{*}(t)}$, by Lemma $2.1(iii)$, there exists $L_0\gg1$ such that $\lambda_1(d, \alpha, \gamma, l, T)\leq \lambda_1(d, \alpha, \min_{[0,T]}\gamma_{*}, l, T)<0$ for all $l\geq L_0$. Consider the following auxiliary problem
\begin{align*}
\left\{\begin{array}{l}
\underline{P}_t=d\underline{P}_{xx}-\alpha \underline{P}_x+\underline{P}(\gamma(x,t)-\underline{P}), \quad 0<x<l,~0\leq t\leq T,\\[5pt]
\underline{P}_{x}(0, t)=\underline{P}(l,t)=0, \quad 0\leq t\leq T,\\[5pt]
\underline{P}(x, 0)=\underline{P}(x, T),\quad 0\leq x \leq l.
\end{array}\right.
 \tag{3.5}
\end{align*}
For such $l$, utilizing \cite{hess91}, $(3.5)$ admits a unique positive solution, denoted by $\underline{P}^l(x,t)$. The comparison principle and the maximum principle yield that $\underline{P}^l$ is increasing in $l$ and
$\underline{P}^l\leq \|\gamma\|_{\infty}:=\|\gamma\|_{L^{\infty}([0,\infty)\times[0,T])}$. Using the regularity theory for parabolic equations, we can show that $\underline{P}^l$ converges to a positive bounded solution of $(3.2)$ as $l\rightarrow \infty$, denoted by $P_*(x,t)$. Let $P(x,t)$ be a positive solution of $(3.2)$.
Since $\underline{P}^l(l,t)=0<P(l,t)$, by the comparison principle, we have $\underline{P}^l(x,t)\leq P(x,t)$ in $[0, l]\times[0,T]$
for any $l\geq L_0$, which implies that $P_*(x,t)\leq P(x,t)$ in $[0, \infty)\times[0,T]$.
So $P_*(x,t)$ is the minimal positive solution of $(3.2)$.

\textbf{Step 2.} some a priori estimates for any positive solution of $(3.2)$
as $x\rightarrow \infty$.

Consider the following problem
\begin{align*}
\left\{\begin{array}{l}
q_t=dq_{xx}-\alpha q_x+q(\gamma_*(t)-q), \quad 0<x<\infty,~0\leq t\leq T,\\[5pt]
q(0, t)=0, \quad 0\leq t\leq T,\\[5pt]
q(x, 0)=q(x, T),\quad 0\leq x <\infty.
\end{array}\right.
 \tag{3.6}
\end{align*}
According to Propositions 2.1 and 2.3 in \cite{dgp13}, $(3.6)$
admits a unique positive solution, denoted by $q(x,t)$, if and only if $0\leq\alpha<\frac{4d}{T}\int_0^T\gamma_{*}(t)dt$. Moreover,
$q_{x}(x,t)>0$, $q(x,t)\rightarrow \underline{q}(t)$ uniformly for $t\in[0,T]$ as $x\rightarrow+\infty$, where
$\underline{q}(t)$ is the unique positive solution of
\begin{align*}
\underline{q}_{t}=\underline{q}(\gamma_*(t)-\underline{q}), ~\mbox{in}~[0, T],\quad \underline{q}(0)=\underline{q}(T).
\end{align*}
For $l\geq L_0$, let us consider the problem
\begin{align*}
\left\{\begin{array}{l}
\mathfrak{q}_t=d\mathfrak{q}_{xx}-\alpha \mathfrak{q}_x+\mathfrak{q}(\gamma_{*}(t)-\mathfrak{q}), \quad 0<x<l,~0\leq t\leq T,\\[5pt]
\mathfrak{q}(0,t)=\mathfrak{q}(l,t)=0,\quad 0\leq t\leq T,\\[5pt]
\mathfrak{q}(x,0)=\mathfrak{q}(x,T),\quad0\leq x\leq l.
\end{array}\right.
 \tag{3.7}
\end{align*}
For such $l$, it follows from \cite{hess91} that $(3.7)$ admits a unique positive solution, denoted by $\mathfrak{q}^{l}(x,t)$. Using the Hopf's boundary lemma, we obtain $\mathfrak{q}_{x}^{l}(0,t)>0$. Then by the comparison principle and strong maximum principle, we obtain that $P\geq\mathfrak{q}^l>0$ in $(0, l)\times[0,T]$ and $\mathfrak{q}^l$ is inceasing in $l$. Similar as in Step 1, we can derive that $\mathfrak{q}^l\rightarrow q$ in $C^{2,1}([0, L]\times[0,T])$ for any $L>0$ as $l\rightarrow \infty$ by using the fact that $q(x,t)$ is the unique positive solution of $(3.6)$. Thus, we have $P(x,t)\geq q(x,t)>0$ in $[0, \infty)\times[0,T]$ and then $\liminf_{x\rightarrow\infty}P(x,t)\geq \underline{q}(t)>0$ uniformly in $[0,T]$.

\textbf{Step 3.} the uniqueness of positive bounded solution of $(3.2)$.

Arguing indirectly, we assume that $(3.2)$ has a positive bounded solution $P(x,t)$ such that $P\not\equiv P_*$. By the strong maximum principle, we easily deduce that $P>P_*$ in $[0,\infty)\times[0,T]$.
Due to the estimate in Step 2 and $P, P_{*}$ is bounded, we can find a constant $k>1$
such that $P\leq kP_*$ in $[0,\infty)\times[0,T]$. Now we turn to a technique introduced by Marcus and V\'{e}ron in \cite{mv98}. Define $Q=P_*-(2k)^{-1}(P-P_*)$.
By direct computations, we have
\begin{align*}
P_*>Q\geq\frac{k+1}{2k}P_*,\quad \frac{2k}{2k+1}Q+\frac{1}{2k+1}P=P_*.
\tag{3.8}
\end{align*}
Noticing that $P(P-\gamma)$ is convex in $P\in(0, \infty)$, we have $P_*^2\leq\frac{2k}{2k+1}Q^2+\frac{1}{2k+1}P^2$ by $(3.8)$. Then direct computation gives $Q_t\geq dQ_{xx}-\alpha Q_x+Q(\gamma(x,t)-Q)$ in $[0,l]\times[0,T], Q_x(0,t)=0, Q(l,t)>0$ in $[0, T]$ and $Q(x, 0)=Q(x,T)$ in $[0, l]$.
This indicates that $Q$ is a super-solution of $(3.5)$. It follows from the comparison principle that $\underline{P}^l\leq Q$ in $[0,l]\times[0,T]$. Since $\underline{P}^l\rightarrow P_*$ in $C^{2,1}([0, L]\times[0, T])$ for any $L>0$ as $l\rightarrow+\infty$.
Hence $P_*\leq Q$ in $[0,\infty)\times[0,T]$. This is a contradiction with $(3.8)$. So $P=P_*$ and the uniqueness is derived.

\textbf{Step 4.} the limit superior of $P$ for large $x$.

By the definition of $\gamma^{\infty}(t)$, we may
assume that $\gamma(x,t)\leq\gamma^{\infty}(t)+\varepsilon
:=\gamma^{\infty}_{\varepsilon}(t)$ for $x\geq L_{0}$.

We consider
\begin{align*}
\left\{\begin{array}{l}
\textsc{p}_t=d\textsc{p}_{xx}-\alpha \textsc{p}_x+\textsc{p}(\gamma_{\varepsilon}^{\infty}(t)-\textsc{p}), \quad L_0<x<l,~0\leq t\leq T,\\[5pt]
\textsc{p}(L_0, t)=\underline{P}^l(L_0,t),~\textsc{p}(l,t)=0, \quad 0\leq t\leq T,\\[5pt]
\textsc{p}(x, 0)=\textsc{p}(x, T),\quad L_0\leq x \leq l,
\end{array}\right.
 \tag{3.9}
\end{align*}
where $\underline{P}^l(x,t)$ is the unique positive solution of $(3.5)$. Clearly, $\underline{P}^l$ is a lower solution of $(3.9)$ and any large constant $M$ is an upper solution of $(3.9)$. Hence it has a positive bounded solution. Since the nonlinear term is concave, we can prove that the positive solution, denoted by $\textsc{p}_{l}^{\varepsilon}(x,t)$, is unique. Therefore,
$\underline{P}^l(x,t)\leq \textsc{p}_{l}^{\varepsilon}(x,t)\leq M$
for any $(x,t)\in [L_{0}, l]\times [0, T]$.
By the regularity theory we have $\textsc{p}_{l}^{\varepsilon}\rightarrow\bar{\textsc{p}}^{\varepsilon}$ locally uniformly in $[L_0, \infty)\times[0,T]$ as $l\rightarrow\infty$, where $\bar{\textsc{p}}^{\varepsilon}$
is a positive bounded solution of
\begin{align*}
\left\{\begin{array}{l}
\bar{\textsc{p}}_t=d\bar{\textsc{p}}_{xx}-\alpha \bar{\textsc{p}}_x+\bar{\textsc{p}}(\gamma_{\varepsilon}^{\infty}(t)-\bar{\textsc{p}}), \quad L_0<x<\infty,~0\leq t\leq T,\\[5pt]
\bar{\textsc{p}}(L_0, t)=P(L_0,t),  \quad 0\leq t\leq T,\\[5pt]
\bar{\textsc{p}}(x, 0)=\bar{\textsc{p}}(x, T),\quad L_0\leq x < \infty.
\end{array}\right.
\tag{3.10}
\end{align*}
Recalling that
$\underline{P}^l(x,t)\rightarrow P(x,t)$ as $l\rightarrow \infty$, we have
$P(x,t)\leq \bar{\textsc{p}}^{\varepsilon}(x,t)$ for $(x,t)\in[L_0, \infty)\times[0, T]$.
We can also show that $\bar{\textsc{p}}^{\varepsilon}(x,t)\rightarrow \mathfrak{p}_{\varepsilon}^{\infty}(t)$ uniformly for $t\in [0, T]$ as $x\rightarrow \infty$,
where $\mathfrak{p}_{\varepsilon}^{\infty}$ is the unique positive solution of $(3.4)$
with $\gamma^{\infty}(t)$ replaced by $\gamma_{\varepsilon}^{\infty}(t)$.
In fact, for any sequence $\{x_{n}\}$ with $x_{n}\rightarrow \infty$ as $n\rightarrow \infty$,
we define $\bar{\textsc{p}}_{n}^{\varepsilon}(x,t)=\bar{\textsc{p}}^{\varepsilon}(x+x_{n},t)$,
then $\bar{\textsc{p}}_{n}^{\varepsilon}$ solves the same equation as $\bar{\textsc{p}}^{\varepsilon}$ but over $(-x_{n}+L_{0}, \infty)\times (0, T)$.
Since $\bar{\textsc{p}}_{n}^{\varepsilon}(x,t)\leq M$, the standard regularity argument allows
us to conclude that we can extract a subsequence of $\{\bar{\textsc{p}}_{n}^{\varepsilon}\}$
(still denoted by $\{\bar{\textsc{p}}_{n}^{\varepsilon}\}$) such that
$\bar{\textsc{p}}_{n}^{\varepsilon}\rightarrow \tilde{\textsc{p}}$ locally in
$C^{2,1}((-\infty, \infty)\times [0, T])$ as $n\rightarrow \infty$
and $\tilde{\textsc{p}}$ is a positive bounded solution of
\begin{align*}
\left\{\begin{array}{l}
w_t=dw_{xx}-\alpha w_x+w(\gamma_{\varepsilon}^{\infty}(t)-w), \quad -\infty<x<\infty,~0<t< T,\\[5pt]
w(x, 0)=w(x, T),\quad -\infty< x < \infty.
\end{array}\right. \tag{3.11}
\end{align*}

Now we show that the positive bounded solution of $(3.11)$ is unique for $\alpha^{2}<4d\min_{t\in[0,T]}\gamma_{\varepsilon}^{\infty}(t)$.
We consider the following boundary-value problem
\begin{align*}
\left\{\begin{array}{l}
-dw_{xx}+\alpha w_x=w(\min_{t\in[0,T]}\gamma_{\varepsilon}^{\infty}(t)-w), \quad -l<x<l,\\[5pt]
w(-l)=w(l)=0.
\end{array}\right. \tag{3.12}
\end{align*}
It is easy to know that $(3.12)$ has a unique positive bounded solution $w^{l}(x)$ for any large $l$ when $\alpha^{2}<4d\min_{t\in[0,T]}\gamma_{\varepsilon}^{\infty}(t)$. Much as before,
we have $w^{l}(x)\rightarrow w(x)$ as $l\rightarrow \infty$ and $w(x)$ is a positive bounded solution of
\begin{align*}
-dw_{xx}+\alpha w_x=w(\min_{t\in[0,T]}\gamma_{\varepsilon}^{\infty}(t)-w), \quad -\infty<x<\infty. \tag{3.13}
\end{align*}
Since $\alpha^{2}<4d\min_{t\in[0,T]}\gamma_{\varepsilon}^{\infty}(t)$, the only positive bounded solution of $(3.13)$ is $w\equiv \min_{t\in[0,T]}\gamma_{\varepsilon}^{\infty}(t)$ (see Remark 4.3 in \cite{bhn05}).
For any positive solution $\hat{\textsc{p}}(x,t)$ of $(3.11)$, by the comparison principle,  $\hat{\textsc{p}}\geq w^{l}$ on $[-l, l]\times [0, T]$.
Letting $l\rightarrow\infty$, we get $\hat{\textsc{p}}\geq \min_{t\in[0,T]}\gamma_{\varepsilon}^{\infty}(t)>0$. Similar to Step 3, we can obtain that the positive bounded solution of $(3.11)$ is unique.

Since the unique positive solution $\mathfrak{p}_{\varepsilon}^{\infty}(t)$ of
\begin{align*}
\mathfrak{p}^{\prime}(t)=\mathfrak{p}(\gamma_{\varepsilon}^{\infty}(t)-\mathfrak{p}),
\quad \mathfrak{p}(0)=\mathfrak{p}(T)
\end{align*}
is bounded and satisfies $(3.11)$, we have $\tilde{\textsc{p}}\equiv \mathfrak{p}_{\varepsilon}^{\infty}$.
Letting $\varepsilon\rightarrow0$, we obtain $\lim_{\varepsilon\rightarrow0}\mathfrak{p}_{\varepsilon}^{\infty}(t)=\mathfrak{p}^{\infty}(t)$ for $t\in[0,T]$.
Thus $\limsup_{x\rightarrow \infty}P(x,t)\leq \mathfrak{p}^{\infty}(t)$ uniformly for $t\in[0,T]$.

For the limit inferior of $P$ for large $x$ can be derived in a similar way, we omit it here.

\textbf{Step 5.} the long time behavior of $p$.

For any $l\geq L_0$, we have $\lambda_1(d, \alpha, \gamma, l, T)<0$. By Theorem 28.1 in \cite{hess91}, we see that the problem
\begin{align*}
\left\{\begin{array}{l}
\underline{\mathcal {P}}_t=d\underline{\mathcal {P}}_{xx}-\alpha \underline{\mathcal {P}}_x+\underline{\mathcal {P}}(\gamma(x,t)-\underline{\mathcal {P}}), \quad 0<x<l,~t\geq0,\\[5pt]
\underline{\mathcal {P}}_x(0, t)=\underline{\mathcal {P}}(l,t)=0,  \quad t\geq0,\\[5pt]
\underline{\mathcal {P}}(x, 0)=p_0(x),\quad 0\leq x \leq l,
\end{array}\right.
\end{align*}
admits a unique positive solution $\underline{\mathcal {P}}^l(x,t)$, which satisfies
$\underline{\mathcal {P}}^l(x, t+nT)\rightarrow \underline{P}^{l}(x,t)$ uniformly for $(x,t)\in[0, l]\times[0, T]$ as $n\rightarrow\infty$, where $\underline{P}^{l}$
is the unique positive solution of $(3.5)$. By the comparison principle, we have $\underline{\mathcal {P}}^l(x, t)\leq p(x,t)$ for $(x,t)\in[0, l]\times[0, \infty)$. Therefore, $\liminf_{n\rightarrow\infty}p(x, t+nT)\geq\underline{P}^{l}(x,t)$ uniformly for $(x,t)\in[0, l]\times[0, T]$. According to Step 1, for such $l$, $\underline{P}^l(x,t)\rightarrow P(x,t)$ in $C^{2,1}([0, L]\times[0,T])$ for any $L>0$
as $l\rightarrow \infty$, where $P(x,t)$ is the unique positive bounded solution of $(3.2)$. Sending $l\rightarrow\infty$, we obtain that $\liminf_{n\rightarrow\infty}p(x, t+nT)\geq P(x,t)$ locally uniformly for $(x,t)\in[0, \infty)\times[0, T]$.

On the other hand, for any $l\geq L_0$, let $\bar{\mathbb{P}}^l(x)$ be the unique positive solution of
\begin{align*}
\left\{\begin{array}{l}
-d\bar{\mathbb{P}}_{xx}+\alpha \bar{\mathbb{P}}_x=\bar{\mathbb{P}}(\|\gamma\|_{\infty}-\bar{\mathbb{P}}), \quad 0<x<l,\\[5pt]
\bar{\mathbb{P}}_x(0)=0, ~\bar{\mathbb{P}}(l)=M,
\end{array}\right.
\end{align*}
where $M$ is a positive constant satisfying $M\geq\|\gamma\|_{\infty}+\|p_0\|_{\infty}$. The existence and uniqueness of $\bar{\mathbb{P}}^l(x)$ can be obtained by the
upper and lower solutions method and the comparison principle.
Thus, $\bar{\mathbb{P}}^l(x)\leq M$ and we can find a constant $k\geq1$ such that $p_{0}(x)\leq k\bar{\mathbb{P}}^l(x)$ for all $0\leq x\leq l$. Note that $p(l,t)\leq kM=k\bar{\mathbb{P}}^l(l)$, by the comparison principle, we have $p(x,t)\leq k\bar{\mathbb{P}}^l(x)$ for $(x,t)\in[0, l]\times[0, \infty)$.

Let $\bar{\mathcal {P}}^l$ be the unique positive solution of
\begin{align*}
\left\{\begin{array}{l}
\bar{\mathcal {P}}_t=d\bar{\mathcal {P}}_{xx}-\alpha \bar{\mathcal {P}}_x+\bar{\mathcal {P}}(\gamma(x,t)-\bar{\mathcal {P}}), \quad 0<x<l,~t\geq0,\\[5pt]
\bar{\mathcal {P}}_x(0, t)=0, \bar{\mathcal {P}}(l,t)=k\bar{\mathbb{P}}^l(l)=kM,  \quad t\geq0,\\[5pt]
\bar{\mathcal {P}}(x, 0)=k\bar{\mathbb{P}}^l(x),\quad 0\leq x \leq l.
\end{array}\right.
\end{align*}
It follows from the comparison principle that $p\leq\bar{\mathcal {P}}^l\leq k\mathbb{P}^l$ for $(x,t)\in[0, l]\times[0, \infty)$. Much as before, we have $\bar{\mathcal {P}}^l(x,t+nT)\rightarrow\bar{P}^l(x,t)$ uniformly for $(x,t)\in[0, l]\times[0, T]$ as $n\rightarrow\infty$, where $\bar{P}^l$ is the unique positive solution of the following problem
\begin{align*}
\left\{\begin{array}{l}
\bar{P}_t=d\bar{P}_{xx}-\alpha \bar{P}_x+\bar{P}(\gamma(x,t)-\bar{P}), \quad 0<x<l,~0\leq t\leq T,\\[5pt]
\bar{P}_x(0, t)=0, \bar{P}(l,t)=kM,  \quad 0\leq t\leq T,\\[5pt]
\bar{P}(x, 0)=\bar{P}(x, T),\quad 0\leq x \leq l.
\end{array}\right.
\end{align*}
Thus we have $\liminf_{n\rightarrow\infty}p(x,t+nT)\leq\bar{P}^l(x,t)$ uniformly for $(x,t)\in[0, l]\times[0, T]$. Note that $\bar{P}^l\leq kM$ and $\bar{P}^l$ is decreasing in $l$. By the regularity theory of parabolic equations, we derive $\bar{P}^l\rightarrow P$ uniformly for $t\in[0, T]$ as $l\rightarrow\infty$. Thus, $\limsup_{n\rightarrow\infty}p(x,t+nT)\leq P(x,t)$ locally uniformly in $[0, \infty)\times[0, T]$. \hfill $\Box$

\subsection{Spreading-vanishing dichotomy of a free boundary problem}

Now we consider the following reaction-diffusion-advection
model with a free boundary in time-periodic environment
\begin{align*}
\left\{\begin{array}{l}
u_t=du_{xx}-\alpha u_x+u(\gamma(x,t)-u), \quad 0<x<h(t),\quad t>0 ,\\[5pt]
u_{x}(0, t)=0,~u(h(t),t)=0,~h'(t)=-\mu u_{x}(h(t), t),\quad t>0,\\[5pt]
h(0)=h_0,~u(x, 0)=u_0(x),\quad 0\leq x< h_0,
\end{array}\right.
 \tag{3.14}
\end{align*}
where $h_0>0, u_0\in C^2([0, h_0])$, and $u_0(x)>0=u_0'(0)=u_0(h_0)$ for $x\in[0, h_0]$,
and assume that $\gamma(x,t)$ satisfies $(H1)$.

Following the arguments of Theorem 3.1
in \cite{dgp13}, we can show that $(3.14)$ admits a unique solution
$(u,h)\in C^{2,1}([0, h(t)]\times [0, \infty))\times C^{1}([0,\infty))$.
Moreover, there exist a constant $M=M(\|\gamma, u_{0}\|_{\infty})>0$ such that
$0<u(x,t)\leq M$ and $0<h^{\prime}(t)\leq \mu M$ for any
$t>0$ and $0<x<h(t)$. Hence $h_{\infty}:=\lim_{t\rightarrow \infty}h(t)$ is well defined. Next, we give the spreading-vanishing dichotomy of $(3.14)$.\\

\noindent\textbf{Proposition 3.2.}
Let $\gamma(x, t)$ be a function satisfying $(H1)$ and
$0<\alpha<2\sqrt{d\min_{[0,T]}\gamma_{*}(t)}$. Then\\
$(i)$ either
\begin{align*}
h_{\infty}=\infty~\mbox{and}~\lim_{n\rightarrow\infty}u(x,t+nT)=P(x,t)~\mbox{locally uniformly for}~(x,t)\in[0,\infty)\times[0,T],
\end{align*}
or
\begin{align*}
h_{\infty}\leq h_*<\infty~\mbox{and}~\limsup_{n\rightarrow\infty}\|u(\cdot,t)\|_{C([0, h(t)])}=0,
\end{align*}
where $P(x,t)$ is the unique positive bounded solution of $(3.2)$ and $h_*$ is the unique positive root of $\lambda_1(d, \alpha, \gamma, \cdot, T)=0$.\\
$(ii)$
when $h_{0}<h_{*}$, then there exists $\mu^{*}$ depending on $(u_{0}, h_{0})$ such that
$h_{\infty}\leq h_{*}$ if $0<\mu\leq\mu^{*}$; and $h_{\infty}=\infty$ if $\mu\geq\mu^{*}$.\\
$(iii)$
when $h_{\infty}=\infty$, we have
\begin{align*}
\frac{1}{T}\int_0^{T}k_0(d, \mu, \gamma_{\infty})(t)dt
\leq\liminf_{t\rightarrow\infty}\frac{h(t)}{t}
\leq\limsup_{t\rightarrow\infty}\frac{h(t)}{t}
\leq\frac{1}{T}\int_0^{T}k_0(d, \mu, \gamma^{\infty})(t)dt+\alpha,
\end{align*}
where $k_0(t)=k_0(d, \mu, \gamma_{\infty})(t)$ is the unique positive continuous $T$-periodic function such that $\mu W_{x}^{k_0}=k_{0}(t)$ in $[0, T]$ and $W^{k_0}(x,t)$ is the unique positive solution of
\begin{align*}
\left\{\begin{array}{l}
W_t=dW_{xx}-k(t)W_x+W(\gamma_{\infty}(t)-W), \quad 0\leq x<\infty,~ 0\leq t\leq T,\\[5pt]
W(0,t)=0,\\[5pt]
W(x, 0)=W(x,T).
\end{array}\right. \tag{3.15}
\end{align*}

\noindent\textbf{Proof.}
We first show $(i)$. Since $h_{\infty}=\infty$,
we know that there exits an integer $m\geq1$ such that $h(t)>l\geq L_0$ for all $t\geq mT$, where $L_{0}$ is defined in the Step 1 of Proposition 3.1. Let $\underline{u}^l$ be the unique positive solution of
\begin{align*}
\left\{\begin{array}{l}
\underline{u}_t=d\underline{u}_{xx}-\alpha \underline{u}_x+\underline{u}(\gamma(x,t)-\underline{u}), \quad 0<x<l,~ t>mT ,\\[5pt]
\underline{u}_{x}(0, t)=0, ~\underline{u}(l,t)=0,\quad t>mT,\\[5pt]
\underline{u}(x, mT)=u(x, mT),\quad 0\leq x\leq l.
\end{array}\right.
\end{align*}
Then, $\underline{u}^l\leq u$ in $[0, l]\times[mT, \infty)$ by using the comparison principle.
It follows from Theorem 28.1 \cite{hess91} that $\underline{u}^l(x, t+nT)\rightarrow \underline{P}^l(x,t)$ uniformly in $[0,l]\times[0, T]$ as $n\rightarrow\infty$, where
$\underline{P}^l(x,t)$ is defined in $(3.5)$. In the Step 3 of Proposition 3.1, we have known $\underline{P}^l\rightarrow P$ uniformly in $[0,L]\times[0, T]$ for
any $L>0$ as $l\rightarrow\infty$. Thus, we obtain that
$\liminf_{n\rightarrow\infty}u(x,t+nT)\geq P(x,t)$ locally uniformly for $(x,t)\in[0,\infty)\times[0,T]$.
On the other hand, similar to Step 5 in Proposition 3.1, we can show that $\limsup_{n\rightarrow\infty}u(x,t+nT)\leq P(x,t)$ locally uniformly for $(x,t)\in[0,\infty)\times[0,T]$. So we get the desired result.

If $h_{\infty}<\infty$, then $h_{\infty}\leq h^*$ and $\lim_{t\rightarrow\infty}\|u(\cdot, t)\|_{C([0, h(t)])}=0$, which can be proved by the similar arguments as Lemma 3.3 in \cite{dgp13} with minor modification, we omit it for brevity.

The proof of $(ii)$ and $(iii)$ is similar to that of Theorems 3.11 and 4.6 in \cite{dgp13},
here we only briefly prove $(iii)$ to give the influence of the advection term.

By the definition of $\gamma_{\infty}$ and $\gamma^{\infty}$, we have
for any small $\varepsilon>0$, there is $L_{*}:=L(\varepsilon)>1$ such that for $x\geq L_{*}$,
\begin{align*}
\gamma(x,t)\leq \gamma_{\varepsilon}^{\infty}(t):= \gamma^{\infty}(t)+\varepsilon,
\quad \gamma(x,t)\geq \gamma_{\infty}^{\varepsilon}(t):= \gamma_{\infty}(t)-\varepsilon.
\end{align*}
By Proposition 3.1, there exists $L^{*}:=L^{*}(\varepsilon)>L_{*}>1$ such that
\begin{align*}
\mathfrak{p}_{\infty}^{\frac{\varepsilon}{2}}(t)
\leq P(x,t)\leq \mathfrak{p}^{\infty}_{\frac{\varepsilon}{2}}(t) \quad
\mbox{for}~(x,t)\in [L^{*}, \infty)\times[0, T],
\end{align*}
where $\mathfrak{p}_{\infty}^{\frac{\varepsilon}{2}}(t)$ and
$\mathfrak{p}^{\infty}_{\frac{\varepsilon}{2}}(t)$ are the positive solutions of $(3.3)$
and $(3.4)$ with $\gamma_{\infty}(t)$ and $\gamma^{\infty}(t)$ replaced by
$\gamma^{\frac{\varepsilon}{2}}_{\infty}(t)$ and $\gamma^{\infty}_{\frac{\varepsilon}{2}}(t)$, respectively.
Since $h_{\infty}=\infty$ and $\lim_{n\rightarrow\infty}u(x,t+nT)=P(x,t)$ locally uniformly in $[0,\infty)\times[0,T]$, there is a positive integer $N=N(L^{*})$ such that with $\mathcal{T}:=NT$,
\begin{align*}
h(\mathcal{T})>3L^{*} \quad \mbox{and}\quad u(2L^{*},t+\mathcal{T})<\mathfrak{p}^{\infty}_{\varepsilon}(t)\quad
\mbox{for all}~ t\geq 0.
\end{align*}
Setting $\tilde{u}(x,t)=u(x+2L^{*}, t+\mathcal{T})$ and
$\tilde{h}(t)=h(t+\mathcal{T})-2L^{*}$. Similar to the proof of Theorem 4.4 (Step 2)
in \cite{dgp13}, there exists $\tilde{\mathcal{T}}=\tilde{\mathcal{T}_{\varepsilon}}=\tilde{N}T>0$
(with an integer $\tilde{N}$) such that
\begin{align*}
\tilde{u}(x,t+\tilde{\mathcal{T}})\leq \mathfrak{p}^{\infty}_{\varepsilon}(t)(1-\varepsilon)^{-1},
\quad \mbox{for}~t\geq \tilde{\mathcal{T}},~x\in[0,\tilde{h}(t)].
\end{align*}

Let $W_{\varepsilon}=W_{\gamma_{\varepsilon}^{\infty}, k^{\varepsilon}}$ be the
unique positive solution of $(3.15)$ with $\gamma^{\infty}(t)=\gamma_{\varepsilon}^{\infty}(t)$
and $k(t)=k^{\varepsilon}(t):=k_{0}(d, \mu, \gamma_{\varepsilon}^{\infty})(t)$.
Since $W_{\varepsilon}(x,t)\rightarrow \mathfrak{p}^{\infty}_{\varepsilon}(t)$ in $[0,T]$ as $x\rightarrow \infty$, there exists $L_{0}^{*}:=L_{0}^{*}(\varepsilon)>2L^{*}$ such that
\begin{align*}
W_{\varepsilon}(x,t)>\mathfrak{p}^{\infty}_{\varepsilon}(t)(1-\varepsilon)
\quad \mbox{for}~(x,t)\in [L_{0}^{*}, \infty)\times[0,T].
\end{align*}
We now define
\begin{align*}
\xi(t)=(1-\varepsilon)^{-2}\int_{0}^{t}(k^{\varepsilon}(s)+\alpha)ds
+L_{0}^{*}+\tilde{h}(\tilde{\mathcal{T}}) \quad\mbox{for}~t\geq 0,\\
w(x,t)=(1-\varepsilon)^{-2}W_{\varepsilon}(\xi(t)-x, t)
\quad \mbox{for}~0\leq x\leq \xi(t), t\geq 0.
\end{align*}
Thus, we have
\begin{align*}
\xi^{\prime}(t)=(1-\varepsilon)^{-2}(k^{\varepsilon}(t)+\alpha),\quad
-\mu w_{x}(\xi(t), t)=\mu(1-\varepsilon)^{-2}(W_{\varepsilon})_{x}(0,t)
=(1-\varepsilon)^{-2}k^{\varepsilon}(t),
\end{align*}
and then
\begin{align*}
\xi^{\prime}(t)\geq-\mu w_{x}(\xi(t), t).
\end{align*}
Direct calculations show that, for $0<x<\xi(t)$ and $t>0$, with $\rho=\xi(t)-x$,
\begin{align*}
w_{t}-dw_{xx}+\alpha w_{x}
&=(1-\varepsilon)^{-2}[(W_{\varepsilon})_{t}+(W_{\varepsilon})_{\rho}\xi^{\prime}(t)
 -d(W_{\varepsilon})_{\rho\rho}-\alpha(W_{\varepsilon})_{\rho}]\\
&\geq(1-\varepsilon)^{-2}[(W_{\varepsilon})_{t}+k^{\varepsilon}(t)(W_{\varepsilon})_{\rho}
 -d(W_{\varepsilon})_{\rho\rho}]\\
&=(1-\varepsilon)^{-2}W_{\varepsilon}(\gamma_{\varepsilon}^{\infty}(t)-W_{\varepsilon})\\
&\geq w(\gamma_{\varepsilon}^{\infty}(t)-w).
\end{align*}
Hence, similar to the proof of Theorem 4.4 (Step 2) in \cite{dgp13}, we have
\begin{align*}
u(x, t+\tilde{\mathcal{T}})\leq w(x,t),~
\tilde{h}(t+\tilde{\mathcal{T}})\leq \xi(t)
\quad \mbox{for}~t\geq 0, ~0\leq x\leq \tilde{h}(t+\tilde{\mathcal{T}}).
\end{align*}
and then
\begin{align*}
\limsup_{t\rightarrow +\infty}\frac{h(t)}{t}
=\limsup_{t\rightarrow +\infty}\frac{\tilde{h}(t-\mathcal{T})+2L^*}{t}
\leq \lim_{t\rightarrow +\infty}\frac{\xi(t-(T+\tilde{T}))+2R^{*}}{t}
\leq \frac{1}{T}\int_{0}^{T}k_{0}(\mu, \gamma^{\infty})(t)dt+\alpha.
\end{align*}

On the other hand, since $\alpha> 0$, we can get
\begin{align*}
\liminf_{t\rightarrow +\infty}\frac{h(t)}{t}
\geq \frac{1}{T}\int_{0}^{T}k_{0}(d, \mu, \gamma_{\infty})(t)dt
\end{align*}
by using the same arguments in the proof of Theorem 4.4 (Step 3) in \cite{dgp13}. \hfill $\Box$

\section{The spreading and vanishing in probelm $(1.1)$}

In this section, we investigate the dynamics of problem $(1.1)$, including the spreading-vanishing quartering, some sufficient conditions for spreading and vanishing, the long-time behavior of solutions and some rough estimates for spreading speed of free boundaries when spreading happens.
To get an entire analysis, we need to add the following assumptions:
\begin{align*}
a_{*}(t)>k^{*}(t)V^{*}(t),\quad b_{*}(t)>h^{*}(t)U^{*}(t),
\tag{H2}
\end{align*}
and
\begin{align*}
0<\alpha_{1}<2\sqrt{d_{1}\min_{[0,T]}[a_{*}(t)-k^{*}(t)V^{*}(t)]},\quad
0<\alpha_{2}<2\sqrt{d_{2}\min_{[0,T]}[b_{*}(t)-h^{*}(t)U^{*}(t)]}, \tag{H3}
\end{align*}
where $U^{*}(t), V^{*}(t)$ are the unique positive solutions of $(3.4)$ with $\gamma^{\infty}$ replaced by $a^{*}(t)$ and $b^{*}(t)$, respectively.

Throughout this section, $(H1)$-$(H3)$ are assumed to hold even if they are not explicitly mentioned.

First, we give the global existence and uniqueness of solutions to the free boundary problem
$(1.1)$. The proof is similar to that of Theorem 1 in \cite{gw15}, we omit the details here. \\

\noindent\textbf{Lemma 4.1.}
Problem $(1.1)$ admits a unique positive solution $(u, v, s_1, s_2)$ which is defined for all $t>0$. Moreover, $s_i\in C^{1+\frac{\nu_0}{2}}([0, \infty))$ $(i=1,2)$, $u\in C^{2,1}(D_1)\cap C^{1+\nu_0, (1+\nu_0)/2}(\bar{D}_1)$, $v\in C^{2,1}(D_2)\cap C^{1+\nu_0, (1+\nu_0)/2}(\bar{D}_2)$, $D_i:=\{(x, t): 0\leq x\leq s_i(t), t>0\}$ $(i=1,2)$, and the following estimates for solutions hold
\begin{align*}
\begin{array}{l}
0<u(x,t)\leq C_1:=\max\{\|a\|_{\infty}, \|u_0\|_{\infty}\},
\quad x\in[0, s_1(t)),~ t\geq0,\\[5pt]
0<v(x,t)\leq C_2:=\max\{\|b\|_{\infty}, \|v_0\|_{\infty}\},
\quad x\in[0, s_2(t)),~ t\geq0,\\[5pt]
0<s_1'(t)\leq C_3\mu_1,\quad 0<s_2'(t)\leq C_4\mu_2,
\quad t>0,
\end{array}
\end{align*}
where $C_3>0$ depending only on $d_{1}, \alpha_{1}$ and $(u_{0}, s_{1}^{0})$;  and $C_4>0$ depending only on $d_{2}, \alpha_{2}$ and $(v_{0}, s_{2}^{0})$.\\

Since $s_i(t)$ are increasing in time, $s_{i,\infty}:=\lim_{t\rightarrow\infty}s_i(t) (i=1,2)$ is well defined. Thus, we can use the same arguments as Lemma 3.1 in \cite{gw15} to obtain the following results.\\

\noindent\textbf{Lemma 4.2.}
Let $(u, v, s_1, s_2)$ be the unique positive solution of $(1.1)$. If $s_{1,\infty}<\infty$ (resp., $s_{2,\infty}<\infty$), then there exists $C>0$ independent of $t$ such that
\begin{align*}
\|u\|_{C^{1+\nu_0, \frac{1+\nu_0}{2}}([0,s_1(t)]\times[T, \infty))}+\|s_1'\|_{C^{\frac{\nu_0}{2}}([T,\infty))}\leq C \\
(\mbox{resp.,}~ \|v\|_{C^{1+\nu_0, \frac{1+\nu_0}{2}}([0,s_2(t)]\times[T, \infty))}+\|s_2'\|_{C^{\frac{\nu_0}{2}}([T,\infty))}\leq C).
\end{align*}
In particular, $\lim_{t\rightarrow\infty}s_1'(t)=0$ (resp., $\lim_{t\rightarrow\infty}s_2'(t)=0$).\\

The spreading-vanishing quartering is a consequence of the following lemmas.\\

\noindent\textbf{Lemma 4.3.}
$(i)$ $\limsup_{n\rightarrow \infty} u(x, t+nT)\leq P_{1}(x,t)$ locally uniformly
in $[0, \infty)\times [0, T]$,
and $\limsup_{n\rightarrow \infty} v(x, t+nT)\leq P_{2}(x,t)$ locally uniformly
in $[0, \infty)\times [0, T]$, where $P_{1}(x, t)$ and $P_{2}(x, t)$ are the unique
solutions of (3.2) with $(d, \alpha, \gamma(x,t))$ replaced by $(d_{1}, \alpha_1, a(x, t))$ and
$(d_{2}, \alpha_2, b(x,t))$, respectively.\\
$(ii)$ $s_{1, \infty}=\infty$ and
$\liminf_{n\rightarrow\infty}u(x,t+nT)\geq Q_1(x,t)$ locally uniformly in $[0,\infty)\times[0,T]$ provided that $s_{1,\infty}>s_1^{*}$, where $s_1^{*}$ is the unique positive root of $\lambda_{1}(d_{1}, \alpha_1, a-kV^{*}, \cdot, T)=0$ and $Q_1(x,t)$ is the unique positive solution
of $(3.2)$ with $(d, \alpha, \gamma(x,t))=(d_{1}, \alpha_1, a(x, t)-k(x,t)V^{*}(t))$.\\
$(iii)$ $s_{2, \infty}=\infty$ and
$\liminf_{n\rightarrow\infty}v(x,t+nT)\geq Q_2(x,t)$ locally uniformly in $[0,\infty)\times[0,T]$ provided that $s_{2,\infty}>s_2^{*}$, where $s_2^{*}$ is the unique positive root of $\lambda_{1}(d_{2}, \alpha_2, b-hU^{*}, \cdot, T)=0$ and $Q_2(x,t)$ is the unique positive solution
of $(3.2)$ with $(d, \alpha, \gamma(x,t))=(d_{2}, \alpha_2, b(x, t)-h(x,t)U^{*}(t))$.\\

\noindent\textbf{Proof.}
$(i)$ In $(3.1)$, we define $(d, \alpha, \gamma(x,t))=(d_{1}, \alpha_1, a(t,x))$ and
\begin{align*}
p_0(x)=
\left\{\begin{array}{l}
u_0(x)\quad 0\leq x\leq s_1^0,\\[5pt]
0\quad x\geq s_1^0,
\end{array}\right.
\end{align*}
and let $p_1(x,t)$ be the corresponding unique positive solution. By the comparison principle,
$u(x,t)\leq p_{1}(x,t)$ for all $0\leq x\leq s_{1}(t)$ and $t>0$. It follows from Proposition 3.1 that $\lim_{n\rightarrow\infty}p_1(x,t+nT)=P_1(x,t)$ locally uniformly in $[0, \infty)\times[0,T]$, where $P_1(x,t)$ is the unique positive bounded solution of $(3.2)$ with $(d, \alpha, \gamma(x,t))=(d_{1}, \alpha_1, a(x,t))$. Thus, we have
$\limsup_{n\rightarrow\infty}u(x,t+nT)\leq P_1(x,t)$ locally uniformly in
$[0, \infty)\times[0,T]$.
Similarly, we have $\limsup_{n\rightarrow\infty}v(x,t+nT)\leq P_2(x,t)$ locally uniformly in $[0, \infty)\times[0,T]$, where $P_2(x,t)$ is the unique positive solution of $(3.2)$ with $(d, \alpha, \gamma(x,t))=(d_{2}, \alpha_2, b(x,t))$.

$(ii)$ We first prove that $s_{1, \infty}=\infty$. Let $\overline{v}(t)$ be the unique solution of the problem
\begin{align*}
\left\{\begin{array}{l}
\overline{v}_{t}=\overline{v}(b^{*}(t)-\overline{v}),~ t>0,\\[5pt]
\overline{v}(0)=\|v_{0}\|_{\infty},
\end{array}\right.
\end{align*}
we have $\lim_{n\rightarrow \infty}\overline{v}(t+nT)=V^{*}(t)$, where $V^{*}(t)$ is the unique positive solutions of $(3.4)$ with $\gamma^{\infty}$ replaced by $b^{*}(t)$. Moreover, $\overline{v}$ satisfies
\begin{align*}
&\overline{v}_{t}=d_{2}\overline{v}_{xx}-\alpha_2\bar{v}_x+\overline{v}(b^{*}(t)-\overline{v})
\geq d_{2}\overline{v}_{xx}-\alpha_2\bar{v}_x+\overline{v}(b(x,t)-\overline{v}),\\
&\overline{v}_{r}(t)=0=v_{r}(0,t),~\overline{v}(t)\geq 0=v(s_{2}(t), t),
~\overline{v}(0)\geq v(x,0).
\end{align*}
Thus, we can apply the comparison principle to deduce
\begin{align*}
v(x, t)\leq \overline{v}(t) \quad\mbox{for}~ 0<x<s_{2}(t), ~t>0.
\end{align*}
As a consequence,
\begin{align*}
\lim_{n\rightarrow \infty}v(x, t+nT)
\leq \lim_{n\rightarrow \infty}\overline{v}(t+nT)=V^{*}(t)
\quad\mbox{uniformly in}~[0,\infty)\times[0, T].
\tag{4.1}
\end{align*}

Since $(H2)$ and $(H3)$ hold, there exists a sufficiently small $\varepsilon_0>0$ such that $a_*(t)>k^*(t)(V^*(t)+\varepsilon)$
and $0<\alpha_{1}<2\sqrt{d_{1}\min_{[0,T]}[a_{*}(t)-k^{*}(t)(V^*(t)+\varepsilon)]}$ for all $0<\varepsilon\leq \varepsilon_0$. Thus, $a_{\varepsilon}(x,t):=a(x,t)-k(x,t)(V^*(t)+\varepsilon)$ satisfies $(H1)$.
By Lemma 2.1, there exists $\bar{s}_{1,\varepsilon}>0$ such that $\bar{s}_{1,\varepsilon}$ is the unique positive root of $\lambda_1(d_1, \alpha_1, a_{\varepsilon},\cdot,T)=0$.
Since $s_{1,\infty}>s_{1}^{*}$, we may assume $s_{1,\infty}>\bar{s}_{1,\varepsilon}$ for all $0<\varepsilon<\varepsilon_{0}$ by Lemma 2.1 $(i)$. By $(4.1)$, there exists an integer $N\gg1$ such that $s_1(NT)>\bar{s}_{1,\varepsilon}$ and $v(x,t+nT)\leq V^*(t)+\varepsilon$ for all $(x,t)\in[0,\infty)\times[0,T]$ and $n\geq N$.
Thus, $(u, s_1)$ is a super-solution of
\begin{align*}
\left\{\begin{array}{l}
\underline{u}_t= d_1\underline{u}_{xx}-\alpha_1\underline{u}_x+\underline{u}(a(x,t)-k(x,t)(V^*(t)+\varepsilon)-\underline{u}),~x\in[0, \underline{s}_1(t)],~t\geq NT, \\[5pt]
\underline{u}_{x}(0, t)=0,~ \underline{u}(\underline{s}_1(t),t)=0, ~ \underline{s}_1'(t)=-\mu_1\underline{u}_{x}(\underline{s}_1(t),t),~ t\geq NT,\\[5pt]
\underline{s}_1(NT):=s_1(NT),~ \underline{u}(x, NT)=\delta u(x, NT),~0\leq x\leq \underline{s}_1(NT),
\end{array}\right.
\tag{4.2}
\end{align*}
where $\delta>0$ is a sufficiently small constant. Since $\underline{s}_1(NT)>\bar{s}_{1,\varepsilon}$, then Lemma 2.2 and Propositions 3.2 yield that $s_{1,\infty}\geq \underline{s}_{1,\infty}=\infty$ and $\liminf_{n\rightarrow\infty}u(x, t+nT)\geq\lim_{n\rightarrow\infty}\underline{u}(x,t+nT)
=Q_{1, \varepsilon}(x,t)$ locally uniformly in $[0,\infty)\times[0,T]$, where
$Q_{1, \varepsilon}(x,t)$ is the unique positive solution
of $(3.2)$ with $(d, \alpha, \gamma(x,t))=(d_{1}, \alpha_1, a(x, t)-k(x,t)(V^{*}(t)+\varepsilon))$. Letting $\varepsilon\rightarrow 0$, we completes the proof of $(ii)$.

Similar to the above proof, we can easily get (iii), here we omit the details. \hfill $\Box$ \\

\noindent\textbf{Lemma 4.4.}
$(i)$ If $s_{1,\infty}\leq s_{1,*}$, then $\lim_{t\rightarrow\infty}\|u(\cdot,t)\|_{C([0,s_1(t)])}=0$.\\
$(ii)$ If $s_{2,\infty}\leq s_{2,*}$, then $\lim_{t\rightarrow\infty}\|v(\cdot,t)\|_{C([0,s_2(t)])}=0$,
where $s_{1,*}$ and $s_{2,*}$ are the unique positive root of $\lambda_1(d_1, \alpha_1, a, \cdot, T)=0$ and
$\lambda_1(d_2, \alpha_2, b, \cdot, T)=0$, respectively.\\

\noindent\textbf{Proof.}
We only prove $(i)$, since $(ii)$ can be proved in a similar way.
Choose $l\in[s_{1,\infty}, s_{1,*}]$. Let $\bar{u}$ be the unique positive solution of
\begin{align*}
\left\{\begin{array}{l}
\bar{u}_t= d_1\bar{u}_{xx}-\alpha_1\bar{u}_x+\bar{u}(a(x,t)-\bar{u}),~x\in(0, l),~t>0, \\[5pt]
\bar{u}_{x}(0, t)=0,~ \bar{u}(l,t)=0, ~ t>0,\\[5pt]
\bar{u}(x, 0)=\max\{\|u_0\|_{\infty}, \|a\|_{\infty}\},~0\leq x\leq l.
\end{array}\right.
\tag{4.3}
\end{align*}
Since $l\in[s_{1,\infty}, s_{1,*}]$, it follows from Lemma 2.1 that $\lambda_1(d_1, \alpha_1, a, l, T)\geq 0$. According to Theorem 28.1 in \cite{hess91}, we have $\lim_{t\rightarrow\infty}\|\bar{u}(\cdot,t)\|_{C([0,l])}=0$. By the comparison principle, we obtain that $0<u\leq \bar{u}$ in $[0,s_1(t)]\times[0,\infty)$. Thus $\lim_{t\rightarrow\infty}\|u(\cdot,t)\|_{C([0,s_1(t)])}=0$. \hfill $\Box$\\

\noindent\textbf{Lemma 4.5.}
$(i)$ Suppose that $s_{1,\infty}<\infty$. If $s_1(t)\leq s_2(t)$ for all large $t$, then
$\lim_{t\rightarrow\infty}\|u(\cdot,t)\|_{C([0,s_1(t)])}=0$.\\
$(ii)$ Suppose that $s_{2,\infty}<\infty$. If $s_2(t)\leq s_1(t)$ for all large $t$, then
$\lim_{t\rightarrow\infty}\|v(\cdot,t)\|_{C([0,s_2(t)])}=0$.\\

\noindent\textbf{Proof.} We only consider $(i)$, since $(ii)$ is similar.
For contradiction, we assume that
$\varepsilon_{0}=\limsup_{t\rightarrow \infty}\|u(\cdot,t)\|_{C([0, s_{1}(t)])}>0,$
then there exists a sequence $(x_{k}, t_{k})\in [0, s_{1}(t_{k})]\times(0, \infty)$ with
$t_{k}\rightarrow \infty$ as $k\rightarrow \infty$, such that $u(x_{k}, t_{k})\geq\frac{\varepsilon_{0}}{2}$ for all $k\in \mathbb{N}$.
Since $0\leq x_{k}<s_{1,\infty}<\infty$, passing to a subsequence if necessary,
one has $x_{k}\rightarrow x_{0}$ as $k\rightarrow \infty$.
We claim $x_{0}\in [0, s_{1,\infty})$. Indeed, if $x_{0}=s_{1,\infty}$, then
$x_{k}-s_{1}(t_{k})\rightarrow 0$ as $k\rightarrow \infty$. Furthermore,
we have
\begin{align*}
\left|\frac{\varepsilon_{0}}{2(x_{k}-s_{1}(t_{k}))}\right|
\leq \left|\frac{u(x_{k}, t_{k})}{x_{k}-s_{1}(t_{k})}\right|
=\left|\frac{u(x_{k}, t_{k})-u(s_{1}(t_{k}), t_{k})}{x_{k}-s_{1}(t_{k})}\right|
=|u_{x}(\bar{x}_{k},t_{k})|\leq C,
\end{align*}
where $\bar{x}_{k}\in(x_{k}, s_{1}(t_{k}))$. It is a contradiction since
$x_{k}-s_{1}(t_{k})\rightarrow 0$ as $k\rightarrow \infty$.

Since $s_{1,\infty}<\infty$, we can use the following transformation
\begin{align*}
y:=\frac{x}{s_{1}(t)},~(U,V)(y,t):=(u,v)(x,t),
~\eta(t):=\frac{s_{2}(t)}{s_{1}(t)},
\end{align*}
then $(U,V)$ satisfies the following problem
\begin{align*}
\left\{\begin{array}{l}
U_t=\frac{d_1}{(s_{1}(t))^{2}}U_{yy}+(\frac{s'_{1}(t)y}{s_{1}(t)}-\frac{\alpha_1}{s_1(t)})U_{y}
+U(a(s_{1}(t)y,t)-U-k(s_{1}(t)y,t)V), \quad 0<y<1,\quad t>0 ,\\[5pt]
V_t=\frac{d_2}{(s_{1}(t))^{2}}V_{yy}+(\frac{s'_{1}(t)y}{s_{1}(t)}-\frac{\alpha_2}{s_1(t)})V_{y}
+V(b(s_{1}(t)y,t)-V-h(s_{1}(t)y,t)U), \quad 0<y<\eta(t),\quad t>0,\\[5pt]
U_{y}(0, t)=V_{y}(0, t)=0, \quad t>0,\\[5pt]
U\equiv 0,\quad y\geq 1,~ t>0;\quad V\equiv 0,\quad y\geq \eta(t),~ t>0,\\[5pt]
\eta(0)=\frac{s_2^0}{s_1^0},~U(y, 0)=u_0(s_{1}^{0}y),~V(y, 0)=v_0(s_{1}^{0}y), \quad 0\leq y< \infty.
\end{array}\right.
\tag{4.4}
\end{align*}
Since $s_{1}(t)\leq s_{2}(t)$ for all large $t$, we have
$\eta(t)\geq 1$ for all large $t$.
Similar to the proof of Lemma 4.2, we obtain
\begin{align*}
\|V\|_{C^{1+\nu_{0}, \frac{1+\nu_{0}}{2}}([0,1]\times[T,\infty))}
\leq C \quad\mbox{for some positive constant C}.
\end{align*}
In view of Lemma 4.2, we derive
\begin{align*}
\|U\|_{C^{1+\nu_{0}, \frac{1+\nu_{0}}{2}}([0,1]\times[T,\infty))}
+\|V\|_{C^{1+\nu_{0}, \frac{1+\nu_{0}}{2}}([0,1]\times[T,\infty))}
\leq C \quad\mbox{for some positive constant C}.
\end{align*}

We now consider
\begin{align*}
u_{k}(y,t):=U(y,t+t_{k}),~ v_{k}(y,t):=V(y,t+t_{k})
\quad\mbox{for all}~y\in [0,1]~\mbox{and}~t\in[0,T].
\end{align*}
Letting $t_{k}=\bar{t}_{k}+\bar{k}T$ with $\bar{t}_{k}\in [0, T)$ and $\bar{k}\in \mathbb{N}$,
passing to a subsequence if necessary, we may assume that $\bar{t}_{k}\rightarrow t_{0}$
as $k\rightarrow \infty$.
Since $s_{1,\infty}<\infty$, we have
$\lim_{k\rightarrow \infty}s'_{1}(t_{k})=0$.
Therefore, we have, up to a subsequence,
\begin{align*}
(u_{k},v_{k})(y,t)\rightarrow (u^{*}, v^{*})(y,t)
\end{align*}
where $u^{*}(\frac{x_{0}}{s_{1,\infty}}, 0)>0$ and
\begin{align*}
\left\{\begin{array}{l}
u^*_t=\frac{d_1}{(s_{1,\infty})^{2}}u^*_{yy}-\frac{\alpha_1}{s_{1, \infty}}u_y^*
+u^*(a(s_{1,\infty}y,t+t_{0})-u^*-k(s_{1,\infty}y,t+t_{0})v^*), \quad 0<y<1,\quad 0<t<T ,\\[5pt]
u^*_{y}(0, t)=0=u^*(1,t), \quad 0<t<T.
\end{array}\right.
\tag{4.5}
\end{align*}
Then the strong maximum principle implies that
$u^*>0$ over $\{(y,t): y\in (0,1), t\in (0,T)\}$.
By Hopf's boundary lemma, there exists $\theta>0$ such that
$u^*_{y}(1,t)\leq -\theta$ for all $t\in(\frac{T}{4}, T)$.
Thus,
\begin{align*}
s'_{1}(t_{k}+\frac{T}{2})
=-\mu_{1}u_{x}(s_{1}(t_{k}+\frac{T}{2}), t_{k}+\frac{T}{2})
=-\mu_{1}\frac{u_{k,y}(1, \frac{T}{2})}{s_{1}(t_{k}+\frac{T}{2})}
\geq\frac{\theta\mu_{1}}{2s_{1,\infty}}
\end{align*}
for all large $k$. This contradicts with Lemma 4.2. Hence, we must have
$\lim_{t\rightarrow\infty}\|u(\cdot,t)\|_{C([0,s_1(t)])}=0$. \hfill $\Box$\\

\noindent\textbf{Lemma 4.6.}
$(i)$ Suppose that $s_{1,\infty}\in(s_{1,*}, s_1^*]$, then $s_1(t)-s_2(t)$ changes sign only finitely many times for large time. Moreover, $s_{2,\infty}=\infty$ and $\lim_{t\rightarrow\infty}\|u(\cdot,t)\|_{C([0,s_1(t)])}=0, \lim_{n\rightarrow\infty}v(x,t+nT)=P_2(x,t)$ locally uniformly in $[0,\infty)\times[0,T]$.\\
$(ii)$ Suppose that $s_{2,\infty}\in(s_{2,*}, s_2^*]$, then $s_1(t)-s_2(t)$ changes sign only finitely many times for large time. Moreover, $s_{1,\infty}=\infty$ and $\lim_{t\rightarrow\infty}\|v(\cdot,t)\|_{C([0,s_2(t)])}=0, \lim_{n\rightarrow\infty}u(x,t+nT)=P_1(x,t)$ locally uniformly in $[0,\infty)\times[0,T]$,
where $P_1$ and $P_2$ are defined in Lemma 4.3.\\

\noindent\textbf{Proof.}
We only deal with $(i)$, since $(ii)$ can be proved similarly.

We first show that $s_{2,\infty}>s_{2,*}$ by using a contradiction argument. If it is not true, then Lemma 4.3 $(ii)$ implies that $\lim_{t\rightarrow\infty}\|v(\cdot,t)\|_{C([0,s_2(t)])}=0$.
Similar to the proof of (ii) in Lemma 4.2, by using $s_{1,\infty}>s_{1,*}$, we have $s_{1, \infty}=\infty$, a contradiction to the assumption $s_{1, \infty}\leq s_1^*$. Thus we have $s_{2,\infty}>s_{2,*}$.

Next, we use a contradiction argument to show that $s_1(t)-s_2(t)$ changes sign only finite many times for large time. Assume that it changes sign infinite many times for large time, then we have $s_{1,\infty}=s_{2,\infty}<\infty$. If we can show that
\begin{align*}
\lim_{t\rightarrow\infty}\|u(\cdot,t)\|_{C([0,s_1(t)])}=0,
\tag{4.6}
\end{align*}
then we can derive that $s_{2,\infty}=\infty$ by using $s_{2,\infty}>s_{2,*}$. Indeed,
since $\lim_{t\rightarrow\infty}\|u(\cdot,t)\|_{C([0,s_1(t)])}=0$
and $s_{2,\infty}>s_{2,*}$, then for any $\varepsilon>0$, there exists an integer $N\gg1$ such that $u(x,t+nT)\leq\varepsilon$ for all $(x,t)\in[0, \infty)\times[0,T], n\geq N$ and $s_2(NT)>\bar{s}_{2}^{\varepsilon}$ due to $\bar{s}_{2}^{\varepsilon}\rightarrow s_{2, *}$ as $\varepsilon\rightarrow0$, where $\bar{s}_{2}^{\varepsilon}$ is
the unique positive root of $\lambda_1(d_1, \alpha_1, b-k\varepsilon,\cdot,T)=0$.
This implies $(v,s_2)$ is a super-solution of
\begin{align*}
\left\{\begin{array}{l}
\underline{v}_t= d_2\underline{v}_{xx}-\alpha_{2}\underline{v}_{x}+\underline{v}(b(x,t)-k(x,t)\varepsilon-\underline{v}),~x\in[0, \underline{s}_2(t)],~t\geq NT, \\[5pt]
\underline{v}_{x}(0, t)=0,~ \underline{v}(\underline{s}_2(t),t)=0, ~ \underline{s}_2'(t)=-\mu_2\underline{v}_{x}(\underline{s}_2(t),t),~ t\geq NT,\\[5pt]
\underline{s}_2(NT):=s_2(NT),~ \underline{v}(x, NT)=\delta v(x, NT),~0\leq x\leq \underline{s}_2(NT),
\end{array}\right.
\tag{4.7}
\end{align*}
where $\delta>0$ is a sufficiently small constant.
Note that $\varepsilon$ is arbitrary, by Proposition 3.2, we have
$\lim_{n\rightarrow\infty}\underline{v}(x,t+nT)=P_2(x,t)$ locally uniformly $(x,t)\in[0,\infty)\times[0,T]$ and $s_{2,\infty}\geq\underline{s}_{2,\infty}=\infty$.
This is a contradiction. Now, we prove $(4.6)$. The proof is similar to that of Lemma 4.5 with minor modification.
For contradiction, we assume that
$\varepsilon_{0}:=\limsup_{t\rightarrow \infty}\|u(\cdot,t)\|_{C([0, s_{1}(t)])}>0,$
then there exists a sequence $(x_{k}, t_{k})\in [0, s_{1}(t_k)]\times(0, \infty)$ with
$t_{k}\rightarrow \infty$ as $k\rightarrow \infty$ such that $u(x_{k}, t_{k})\geq\frac{\varepsilon_{0}}{2}$ for all $k\in \mathbb{N}$.
Since $0\leq x_{k}<s_{1,\infty}<\infty$, passing to a subsequence if necessary,
one has $x_{k}\rightarrow x_{0}$ as $k\rightarrow \infty$.
Similar to the proof of Lemma 4.5, we can show $x_{0}\in [0, s_{1,\infty})$.

Since $s_{1,\infty}<\infty$, we can use the following transformation
\begin{align*}
y:=\frac{x}{s_{1}(t)},~(U,V)(y,t):=(u,v)(x,t),
~\eta(t):=\frac{s_{2}(t)}{s_{1}(t)},
\end{align*}
then $(U,V)$ satisfies $(4.4)$.

Define $\beta_{k}:=\min\{1, \min_{t\in[t_{k}, t_{k}+T]}\eta(t)\}$, we know
$\lim_{k\rightarrow\infty}\beta_{k}=1$ since $s_{1,\infty}=s_{2,\infty}<\infty$.
We now consider
\begin{align*}
u_{k}(y,t):=U(y,t+t_{k}),~ v_{k}(y,t):=V(y,t+t_{k})
\quad\mbox{for all}~y\in [0,\beta_{k}]~\mbox{and}~t\in[0,T].
\end{align*}
By Lemma 4.2, we know that
\begin{align*}
\|u_{k}\|_{C^{1+\nu_{0}, \frac{1+\nu_{0}}{2}}([0,\beta_{k}]\times[0,T])}
+\|v_{k}\|_{C^{1+\nu_{0}, \frac{1+\nu_{0}}{2}}([0,\beta_{k}]\times[0,T])}\leq C
\end{align*}
for some positive constant $C$ independent of $k$.
Letting $t_{k}=\bar{t}_{k}+\bar{k}T$ with $\bar{t}_{k}\in [0, T)$ and $\bar{k}\in \mathbb{N}$,
passing to a subsequence if necessary, we may assume that $\bar{t}_{k}\rightarrow t_{0}$
as $k\rightarrow \infty$.
Since $s_{1,\infty}<\infty$, we have
$\lim_{k\rightarrow \infty}s'_{1}(t_{k})=0$.
Therefore, we have, up to a subsequence,
\begin{align*}
(u_{k},v_{k})(y,t)\rightarrow (u^{*}, v^{*})(y,t),
\end{align*}
where $u^{*}(\frac{x_{0}}{s_{1,\infty}}, 0)>0$ and satisfies $(4.5)$.
As in the proof of Lemma 4.5, using the strong maximum principle
and Hopf's boundary lemma, we can derive
$s'_{1}(t_{k}+\frac{T}{2})\geq \delta$
for some $\delta>0$ and for all large $k$. This contradicts with Lemma 4.2. Hence, we must have
$\lim_{t\rightarrow\infty}\|u(\cdot,t)\|_{C([0,s_1(t)])}=0$, and then
$s_{1}(t)-s_{2}(t)$ changes sign only finitely many times.

From the above analysis, we see that either $s_{1}(t)\leq s_{2}(t)$ for all large $t$
or $s_{1}(t)\geq s_{2}(t)$ for all large $t$. However, the latter case cannot happen.
Otherwise, we have $s_{2,\infty}<\infty$, and then
$\lim_{t\rightarrow\infty}\|v(\cdot, t)\|_{C([0, s_{2}(t)])}=0$
by Lemma 4.5 $(ii)$. Similar to (4.7), we can construct a sub-solution of $u$ and get $s_{1,\infty}=\infty$, which is a contradiction.
Thus, $s_{1}(t)\leq s_{2}(t)$ for all large $t$. Consequently, Lemma 4.5 $(i)$
implies that $\lim_{t\rightarrow\infty}\|u(\cdot,t)\|_{C([0,s_1(t)])}=0$, and then
$s_{2,\infty}=\infty$.

Finally, we give the long time behavior of $v(x,t)$. By Lemma 4.3 $(i)$, we have
\begin{align*}
\limsup_{n\rightarrow\infty}v(x,t+nT)\leq P_2(x,t)~\mbox{locally uniformly in}~[0, \infty)\times[0,T].
 \tag{4.8}
\end{align*}
Since $\lim_{t\rightarrow\infty}\|u(\cdot,t)\|_{C([0,s_1(t)])}=0$ and
$s_{2,\infty}=\infty$, then by the same arguments as above, we can construct a sub-solution $(\underline{v}, \underline{s}_2)$ given by $(4.7)$. Thus,
it follow from Lemma 2.2 that
\begin{align*}
\liminf_{n\rightarrow\infty}v(x,t+nT)\geq\lim_{n\rightarrow\infty}\underline{v}(x,t+nT)= P_2(x,t)~\mbox{locally uniformly in}~[0, \infty)\times[0,T].
 \tag{4.9}
\end{align*}
Combining $(4.8)$ with $(4.9)$, we have $\lim_{n\rightarrow\infty}v(x,t+nT)=P_2(x,t)$
locally uniformly in $[0, \infty)\times[0,T]$. \hfill $\Box$\\

Combining Lemmas 4.3-4.6, we immediately obtain the following spreading-vanishing quartering.\\

\noindent\textbf{Theorem 4.1.}
The dynamics of $(1.1)$ can be classified into four cases:\\
$(i)$ both two species vanish eventually if $s_{1, \infty}\leq s_{1, *}$ and
$s_{2, \infty}\leq s_{2, *}$.\\
$(ii)$ $u$ spreads successfully and $v$ vanishes eventually if $s_{2, \infty}\leq s_2^*$.\\
$(iii)$ $u$ vanishes eventually and $v$ spreads successfully if $s_{1, \infty}\leq s_1^*$.\\
$(iv)$ both two species spread successfully.\\

\noindent\textbf{Proof.}
Let $(u, v, s_1, s_2)$ be the unique bounded positive solution of $(1.1)$, and we shall divide our proof into three steps.

\textbf{Step 1.} the case $s_{2,\infty}\leq s_{2,*}$.

It follows from Lemma 4.4 $(ii)$ that $v$ vanishes eventually. If $s_{1, \infty}\leq s_{1,*}$, then $u$ vanishes eventually by using Lemma 4.4 $(i)$. Thus $(i)$ is proved. Since $\lim_{t\rightarrow\infty}\|v(\cdot, t)\|_{C([0, s_2(t)])}=0$, then Lemma 4.6 $(i)$ cannot happen provided that $s_{1,\infty}\in (s_{1,*}, s_{1}^{*}]$. This implies that $s_{1, \infty}> s_1^*$. By Lemma 4.3 $(i)$, we have $s_{1,\infty}=\infty$, which means $u$ spreads successfully. Thus the solution satisfies $(ii)$.

\textbf{Step 2.} the case $s_{2,*}<s_{2, \infty}\leq s_2^*$.

By Lemma 4.6 $(ii)$, we obtain that the solution satisfies $(ii)$.

\textbf{Step 3.} the case $s_{2, \infty}>s_2^*$.

Using Lemma 4.3 $(iii)$, we derive that $s_{2, \infty}=\infty$, which means $v$ spread successfully. If $s_{1, \infty}> s_1^*$, then Lemma 4.3 $(i)$ shows that $u$ spreads successfully. Then the solution satisfies $(iv)$. Since $v$ spreads successfully, it follows from Lemma 4.6 $(i)$ and Lemma 4.4 $(i)$ that $\lim_{t\rightarrow\infty}\|u(\cdot,t)\|_{C([0,s_1(t)])}=0$ provide that
$s_{1,\infty}\leq s_1^*$. This suggests that $u$
vanishes eventually. So the solution satisfies $(iii)$. \hfill $\Box$\\

To discuss the long-time behavior of solutions, we need to study the following T-periodic boundary-value problem of the diffusion-advection competition model in the half line
\begin{align*}
\left\{\begin{array}{l}
U_t=d_1U_{xx}-\alpha_1U_x+U(a(x,t)-U-k(x,t)V), \quad 0<x<\infty,~0<t<T,\\[5pt]
V_t=d_2V_{xx}-\alpha_2V_x+V(b(x,t)-V-h(x,t)U), \quad 0<x<\infty,~0<t<T,\\[5pt]
U_{x}(0, t)=V_{x}(0, t)=0, \quad 0<t<T,\\[5pt]
U(x, 0)=U(x, T),~V(x, 0)=V(x, T),\quad 0\leq x <\infty,
\end{array}\right.
 \tag{4.10}
\end{align*}
where $a(x,t)$, $b(x,t)$, $k(x,t)$ and $h(x,t)$ are functions satisfying $(H1)$.
The following lemma is essentially parallel to Theorem 2.1 in \cite{w144}. Here we omit the details of proof.\\

\noindent\textbf{Lemma 4.7.}
If the assumptions $(H1)$-$(H3)$ hold,
then there exist four positive $T$-periodic functions $U_{1}$, $U_{2}$, $V_{1}$,
$V_{2}\in (C^{2+\nu_{0}, 1+\frac{\nu_{0}}{2}}\cap L^{\infty})([0,\infty)\times[0, T])$,
such that both $(U_{1}, V_{2})$ and $(U_{2}, V_{1})$ are positive bounded solutions of $(4.10)$.
Moreover, any positive bounded solution $(U, V)$ of $(4.10)$ satisfies
\begin{align*}
U_{1}(x, t)\leq  U(x, t)\leq U_{2}(x, t),~
V_{1}(x, t)\leq V(x, t)\leq V_{2}(x, t), \quad\mbox{uniformly~in}~[0,L]\times[0, T]
\end{align*}
for any $L>0$.
Moreover, any positive bounded solution $(U, V)$ of $(4.10)$ satisfies
\begin{align*}
U_{*}(t)\leq \liminf_{x\rightarrow\infty}U(x, t)
\leq \limsup_{x\rightarrow\infty}U(x, t)\leq U^{*}(t), \quad
V_{*}(t)\leq \liminf_{x\rightarrow\infty}V(x, t)
\leq \limsup_{x\rightarrow\infty}V(x, t)\leq V^{*}(t),
\end{align*}
where $U^{*}(t), V^{*}(t)$ are the unique positive solutions of $(3.4)$ with $\gamma^{\infty}$ replaced by $a^{*}(t)$ and $b^{*}(t)$, respectively,
and $U_{*}(t), V_{*}(t)$ are the unique positive solutions of
$(3.3)$ with $\gamma_{\infty}$ replaced by $a_{*}(t)-k^*(t)V^{*}(t)$ and $b_{*}(t)-h^*(t)U^{*}(t)$, respectively.\\

Now, we establish the long-time behavior of solutions when spreading occurs.\\

\noindent\textbf{Theorem 4.2.}
If $s_{1, \infty}=s_{2,\infty}=\infty$, then for each $l>0$
\begin{align*}
\begin{array}{l}
U_1(x,t)\leq\liminf_{n\rightarrow\infty}u(x, t+nT)\leq\limsup_{n\rightarrow\infty}u(x, t+nT)\leq U_2(x,t),\\[5pt]
V_1(x,t)\leq\liminf_{n\rightarrow\infty}v(x, t+nT)\leq\limsup_{n\rightarrow\infty}v(x, t+nT)\leq V_2(x,t)
\end{array}
\end{align*}
uniformly in $[0,l]\times[0,T]$, where $U_1, U_2, V_1,V_2$ are given in Lemma 4.7.\\

\noindent\textbf{Proof.}
The proof is a simple modification of that for Theorem 3.2 in \cite{w144}. So we omit it here. \hfill $\Box$\\

The following lemma estimates the asymptotic spreading speed of the free boundaries $s_{i}(t)$ $(i=1,2)$.\\

\noindent\textbf{Lemma 4.8.}
$(i)$ Assume that $s_{1,\infty}=\infty$, then
\begin{align*}
c_{1,*}\leq\liminf_{t\rightarrow\infty}\frac{s_1(t)}{t}
\leq\limsup_{t\rightarrow\infty}\frac{s_1(t)}{t}\leq c_{1}^{*},
\tag{4.11}
\end{align*}
where $c_{1,*}=\frac{1}{T}\int_{0}^{T}k_{0}(d_{1}, \mu_{1}, a_{\infty}-k^{\infty}V^{*})(t)dt$ and $c_{1}^{*}=\frac{1}{T}\int_{0}^{T}k_{0}(d_{1}, \mu_{1}, a^{\infty})(t)dt+\alpha_{1}$.\\
$(ii)$ Assume that $s_{2,\infty}=\infty$, then
\begin{align*}
c_{2,*}\leq\liminf_{t\rightarrow\infty}\frac{s_2(t)}{t}
\leq\limsup_{t\rightarrow\infty}\frac{s_2(t)}{t}\leq c_{2}^{*},
 \tag{4.12}
\end{align*}
where $c_{2,*}=\frac{1}{T}\int_{0}^{T}k_{0}(d_{2}, \mu_{2}, b_{\infty}-h^{\infty}U^{*})(t)dt$ and $c_{2}^{*}=\frac{1}{T}\int_{0}^{T}k_{0}(d_{2}, \mu_{2}, b^{\infty})(t)dt+\alpha_{2}$.\\

\noindent\textbf{Proof.}
We only deal with $(i)$, since $(ii)$ can be proved in a similar way.

It is easily to check that $(u, s_1)$ forms a sub-solution of
\begin{align*}
\left\{\begin{array}{l}
\bar{u}_t=d_{1}\bar{u}_{xx}-\alpha_{1}\bar{u}_x
+\bar{u}(a(x,t)-\bar{u}), \quad 0<x<\bar{s}_{1}(t),\quad t>0 ,\\[5pt]
\bar{u}_{x}(0, t)=0, \bar{u}(\bar{s}_{1}(t),t)=0, \bar{s}_{1}'(t)=-\mu_{1} \bar{u}_{x}(\bar{s}_{1}(t), t),\quad t>0,\\[5pt]
\bar{s}_{1}(0)=s_{1}^0,~\bar{u}(x, 0)=u_0(x),\quad 0\leq x<\bar{s}_{1}^0,
\end{array}\right. \tag{4.13}
\end{align*}
and we denote the corresponded solution by $(\bar{u}, \bar{s}_1)$. It follows from Lemma 2.2 that $s_1(t)\leq\bar{s}_1(t)$ for all $t>0$, which implies that $\bar{s}_{1,\infty}=\infty$. Thus by Proposition 3.2 $(iii)$, we have $\limsup_{t\rightarrow\infty}\frac{s_1(t)}{t}
\leq \lim_{t\rightarrow\infty}\frac{\bar{s}_1(t)}{t}
\leq \frac{1}{T}\int_{0}^{T}k_{0}(d_{1}, \mu_{1}, a^{\infty})(t)dt+\alpha_{1}:=c_{1}^{*}$.

Next we give the lower bound estimate in $(4.11)$. Since $s_{1,\infty}=\infty$ and $(4.1)$, for any small $\varepsilon>0$ we can choose an integer $N\gg 1$ such that $s_1(NT)>\bar{s}_{1,\varepsilon}$ and $v(x,t+nT)\leq V^{*}(t)+\varepsilon$ for all $(x,t)\in[0,\infty)\times[0,T]$ and $n\geq N$, where $\bar{s}_{1,\varepsilon}$ is defined in the proof of Lemma 4.3 $(ii)$.
Then it is easily to check that $(u,s_1)$ forms a super-solution of $(4.2)$,
and denote the unique positive solution by $(\underline{u}, \underline{s}_{1})$.
Similar to the proof of Lemma 4.3 $(ii)$, we can obtain $\underline{s}_{1,\infty}=\infty$. It follows from Lemma 2.2 and Propositions 3.2 that
$\lim_{t\rightarrow\infty}\frac{s_1(t)}{t}\geq
\lim_{t\rightarrow\infty}\frac{\underline{s}_1(t)}{t}\geq
\frac{1}{T}\int_{0}^{T}k_{0}(d_{1}, \mu_{1}, a_{\infty}-k^{\infty}(V^{*}+\varepsilon))(t)dt
:=c_{1,*}^{\varepsilon}$.
Letting $\varepsilon\rightarrow 0$, we have
$\lim_{t\rightarrow\infty}\frac{s_1(t)}{t}\geq c_{1,*}$. \hfill $\Box$\\

Next, we are in a position to show the existence of a minimal habitat size for spreading $s_{1,min}$ (resp., $s_{2,min}$) for species $u$ (resp., $v$), which means it is the minimal value such that $s_1^0\geq s_{1,min}$ (resp., $s_2^0\geq s_{2,min}$) guarantees the spreading of $u$ (resp., $v$), regardless of $s_2^0$ (resp., $s_1^0$), $u_0, v_0$ and $\mu_i$ $(i=1,2)$, but it can vanish eventually if $s_1^0<s_{1, min}$ (resp., $s_2^0<s_{2,min}$). Also, an upper/lower estimate of $s_{i,min}$ $(i=1,2)$ can be given.\\

\noindent\textbf{Theorem 4.3.}
$(i)$ There exists minimal habitat size for spreading
\begin{align*}
s_{1, min}:=\min \{\hat{s}_1>0~| ~u~\mbox{always spreads successfully if}~ s_1^0=\hat{s}_1\}
\end{align*}
such that $u$ spreads successfully, regardless of $u_0, v_0, s_2^0$ and $\mu_i$ $(i=1,2)$ if and only if $s_1^0\geq s_{1, min}$. Furthermore, $s_{1, **}\leq s_{1, min}\leq s_1^*$, where $s_{1,**}$ satisfies $\lambda_1(d_1, \alpha_1, a-kV_{1},\cdot,T)=0$.\\
$(ii)$ There exists minimal habitat size for spreading
\begin{align*}
s_{2, min}:=\min \{\hat{s}_2>0~| ~v~\mbox{always spreads successfully if}~ s_2^0=\hat{s}_2\}
\end{align*}
such that $v$ spreads successfully, regardless of $u_0, v_0, s_1^0$ and $\mu_i$ $(i=1,2)$ if and only if $s_2^0\geq s_{2, min}$. Furthermore, $s_{2, **}\leq s_{2, min}\leq s_2^*$, where $s_{2,**}$ satisfies $\lambda_1(d_2, \alpha_2, b-hU_{1},\cdot,T)=0$.\\

\noindent\textbf{Proof.}
We only deal with $(i)$, since $(ii)$ can be proved in a similar way.

\textbf{Step 1.} we shall show that $u$ spreads successfully if and only if $s_1^0\geq s_{1, min}$.

We define $A:=\{\hat{s}_1>0 ~|~ u~\mbox{always spreads successfully if}~s_1^0=\hat{s}_1\}$. Note that $A\neq\varnothing$ since $s_1^*\in A$. Thus, $s_{1, min}:=\inf A$ is well defined.
If $\tilde{s}_1\in A$ for some $\tilde{s}_1>0$, then by Lemma 2.2 we see that $s_1\in A$ for all $s_1>\tilde{s}_1$. Thus, $u$ always spreads successfully for all $s_1^0>s_{1, min}$. Similar to the proof of Theorem 2 in \cite{chw15}, we can show that $s_{1, min}\in A$. Thus,
$u$ spreads successfully, regardless of $u_0, v_0, s_2^0$ and $\mu_i$ $(i=1,2)$ if and only if $s_1^0\geq s_{1, min}$.

\textbf{Step 2.} we will prove that $s_{1,**}\leq s_{1, min}\leq s_1^*$.

Since $s_1^*\in A$, so $s_{1, min}\leq s_1^*$. For the lower bound, we use a contradiction argument. Assume that $s_{1,min}<s_{1,**}$, since $s_{1,min}\in A$ and by choosing $s_1^0=s_{1,min}$, we can obtain that $s_{1,\infty}^{\mu_1, \mu_2}=\infty$.
Now we fix initial data and all parameter of $(1.1)$ except $\mu_1$ and $\mu_2$.

Taking any $\hat{\mu}_1>0$ and $\mu_2>0$, also choosing $s_2^0>s_2^*$ such that $v^{\hat{\mu}_1, \mu_2}$ spreads successfully. Moreover, Theorem 4.2 yields that there exists an integer $N\gg1$ such that
$v^{\hat{\mu}_1, \mu_2}(x,t+nT)\geq V_{1}(x,t)-\varepsilon$ for all $(x,t)\in[0, s_1^*]\times[0,T]$ and $n\geq N$, where $\varepsilon>0$ can be chosen smaller if necessary. By Corollary 2.1, if $0<\mu_1\leq\hat{\mu}_1$, then
\begin{align*}
v^{\mu_1, \mu_2}(x,t+nT)\geq V_{1}(x,t)-\varepsilon \quad\mbox{for}~ (x,t)\in[0, s_1^*]\times[0,T],~ n\geq N. \tag{4.14}
\end{align*}
Using $s_1^0=s_{1,min}$ and $s_{1, min}<s_{1,**}$, we can obtain that $s_1^0<s_{1,**}^\varepsilon$, where $s_{1,**}^\varepsilon$ is the unique positive root of
$\lambda_1(d_1, \alpha_1, a-k(V_{1}-\varepsilon),\cdot,T)=0$.
For such fixed $N$, by Lemma 4.1 ($0<s_{1}^{\prime}(t)\leq C_{3}\mu_{1}$) and Corollary 2.1, there exists $\check{\mu}_1\in(0, \hat{\mu}_1]$ small enough, such that $s_1^{\mu_1, \mu_2}(NT)<s_{1,**}^\varepsilon$ for all $\mu_1\in(0, \check{\mu}_1]$.

Next, we apply comparison principle to derive that $u^{\mu_1, \mu_2}$ vanishes eventually provided that $\mu_1$ is small enough.
For each $\mu_1\leq\check{\mu}_1$, due to (4.14), we have
\begin{align*}
u_{t}^{\mu_{1},\mu_{2}}\leq d_{1}u_{xx}^{\mu_{1},\mu_{2}}-\alpha_{1}u_{x}^{\mu_{1},\mu_{2}}
+u^{\mu_{1},\mu_{2}}(a(x,t)-k(x,t)(V_{1}(x,t)-\varepsilon)-u^{\mu_{1},\mu_{2}})
\end{align*}
for all $x\in [0, \min\{s_{1}^{\mu_{1},\mu_{2}}(t), s_{1}^{*}\}]$ and $t\geq NT$.
Let $(\bar{u}, \bar{s}_{1})$ be the solution of
\begin{align*}
\left\{\begin{array}{l}
\bar{u}_t=d_{1}\bar{u}_{xx}-\alpha_{1}\bar{u}_x
+\bar{u}(a(x,t)-k(x,t)(V_{1}(x,t)-\varepsilon)-\bar{u}),
\quad 0<x<\bar{s}_{1}(t),\quad t>NT,\\[5pt]
\bar{u}_{x}(0, t)=0, ~\bar{u}(\bar{s}_{1}(t),t)=0, ~\bar{s}_{1}'(t)=-\mu_{1} \bar{u}_{x}(\bar{s}_{1}(t), t),\quad t>NT,\\[5pt]
\bar{s}_{1}(NT)=\bar{s}_{1}^{NT},~\bar{u}(x, NT)=\bar{u}^{NT}(x),\quad 0\leq x\leq\bar{s}_{1}^{NT},
\end{array}\right.
\end{align*}
where initial data $(\bar{s}_{1}^{NT}, \bar{u}^{NT}(x))$ satisfies $\bar{s}_{1}^{NT}\in (s_1^{\check{\mu}_1, \mu_2}(NT), s_{1,**}^\varepsilon)$ and
$\bar{u}^{NT}(x)>u^{\check{\mu}_1, \mu_2}(x, NT)$ for $0< x<s_1^{\check{\mu}_1, \mu_2}(NT)$.
Since $\bar{s}_{1}^{NT}<s_{1,**}^\varepsilon$, by Proposition 3.2 $(ii)$, there exists $\mu_{1}^{*}$ such that $\bar{s}_{1,\infty}\leq s_{1,**}^\varepsilon$ and
$\lim_{t\rightarrow \infty}\|\bar{u}(\cdot, t)\|_{C([0, \bar{s}_{1}(t)]}=0$
if $\mu_{1}\leq \mu_{1}^{*}$. By Corollary 2.1, we have
$s_{1}^{\mu_{1},\mu_{2}}(NT)\leq s_{1}^{\check{\mu}_{1},\mu_{2}}(NT)<\bar{s}_{1}^{NT}$
and $u^{\mu_{1},\mu_{2}}(x, NT)\leq u^{\check{\mu}_{1},\mu_{2}}(x, NT)<\bar{u}^{NT}(x)$
for any $\mu_{1}\in(0, \min\{\check{\mu}_{1}, \mu_{1}^{*}\}]$
and $x\in [0, s_{1}^{\mu_{1},\mu_{2}}(NT)]$.
By Lemma 2.2, we can derive that $s_{1, \infty}^{\mu_1, \mu_2}<\bar{s}_{1,\infty}\leq s_{1,**}^\varepsilon$. This contradicts to $s_{1, \infty}^{\mu_1, \mu_2}=\infty$.
Hence, $s_{1,**}\leq s_{1, min}\leq s_1^*$, which completes the proof of $(i)$. \hfill $\Box$\\

Finally, we aim to use parameter $\mu_i$ $(i=1,2)$ (assuming the others are fixed) to derive sharp criteria for spreading and vanishing. Hereafter, we first give two lemmas.\\

\noindent\textbf{Lemma 4.9.}
$(i)$ If $0<s_1^0<s_{1,*}$, then there exists $\underline{\mu}_1>0$ depending on $u_0$ and $s_1^0$ such that $s_{1,\infty}<\infty$ provided that $0<\mu_1<\underline{\mu}_1$.\\
$(ii)$ If $0<s_1^0<s_{1,min}$, then there exists $\bar{\mu}_1>0$ depending on $u_0, v_0, d_1$ and $s_1^0$ such that $s_{1,\infty}=\infty$ provided that $\mu_1>\bar{\mu}_1$.\\
$(iii)$ If $0<s_2^0<s_{2,*}$, then there exists $\underline{\mu}_2>0$ depending on $v_0$ and $s_2^0$ such that $s_{2,\infty}<\infty$ provided that $0<\mu_2<\underline{\mu}_2$.\\
$(iv)$ If $0<s_2^0<s_{2,min}$, then there exists $\bar{\mu}_2>0$ depending on $u_0, v_0, d_2$ and $s_2^0$ such that $s_{2,\infty}=\infty$ provided that $\mu_2>\bar{\mu}_2$.\\

\noindent\textbf{Proof.}
We only deal with $(i)$ and $(ii)$, since $(iii)$ and $(iv)$ can be proved in a similar way.

$(i)$
In the proof of Lemma 4.8, we know that $(u, s_{1})$ is a sub-solution of (4.13),
and then $u(x,t)\leq\bar{u}(x,t), s_1(t)\leq\bar{s}_1(t)$ for all $(x,t)\in[0, s_{1}(t))\times[0, \infty)$. By Proposition 3.2 (ii), there exists $\underline{\mu}_1>0$ depending on $u_0$ and $s_1^0$ such that $\bar{s}_{1,\infty}<\infty$ provided that $0<\mu_1\leq\underline{\mu}_1$. This implies that $s_{1,\infty}\leq \bar{s}_{1,\infty}<\infty$ provided that $0<\mu_1\leq\underline{\mu}_{1}$.

$(ii)$
By $(4.1)$, we have $\lim_{n\rightarrow \infty} v(x, t+nT)\leq V^{*}(t)$ uniformly in
$[0, \infty)\times [0, T]$, then for any small $\varepsilon>0$, there exists an integer
$N\gg 1$ such that $v(x, t+nT)\leq V^{*}(t)+\varepsilon$ for $n\geq N$ and
$(x, t)\in [0, \infty)\times[0, T]$. Thus, $(u, s_{1})$ satisfies
\begin{align*}
\left\{\begin{array}{l}
u_t\geq d_1u_{xx}-\alpha_1u_x+u(a(x,t)-k(x,t)(V^*(t)+\varepsilon)-u),~x\in[0, s_1(t)),~t> NT, \\[5pt]
u_{x}(0, t)=0,~ u(s_1(t),t)=0, ~ s_1'(t)=-\mu_1u_{x}(s_1(t),t),~ t> NT,\\[5pt]
u(x, NT)>0,~0\leq x< s_1(NT).
\end{array}\right.
\end{align*}
Note that $u(x, NT)$ depends on $\mu_{1}$, but $u(x, NT)\geq z(x, NT)$, where
$z(x,t)$ and $w(x,t)$ are determined by
\begin{align*}
\left\{\begin{array}{l}
z_t=d_1z_{xx}-\alpha_1 z_x+z(a(x,t)-z-k(x,t)w), \quad 0<x<s_1^{0},\quad t>0 ,\\[5pt]
w_t=d_2w_{xx}-\alpha_2 w_x+w(b(x,t)-w-h(x,t)z), \quad 0<x<s_1^{0},\quad t>0 ,\\[5pt]
z_{x}(0, t)=w_{x}(0, t)=0, \quad t>0,\\[5pt]
z(s_1^{0}, t)=0, ~ w(s_1^{0}, t)=\max\{\|b\|_{\infty}, \|v_{0}\|_{\infty}\}, \quad t>0,\\[5pt]
z(x, 0)=u_0(x),~w(x, 0)=\max\{\|b\|_{\infty}, \|v_{0}\|_{\infty}\}, \quad 0\leq x\leq s_1^{0}.
\end{array}\right.
\end{align*}
Clearly, $z(x, NT)$ is independent of $\mu_{1}$. Now it is easy to see that $(u, s_{1})$ is an upper-solution to the problem
\begin{align*}
\left\{\begin{array}{l}
\underline{u}_t= d_1\underline{u}_{xx}-\alpha_1\underline{u}_x+u(a(x,t)-k(x,t)(V^*(t)+\varepsilon)-\underline{u}),~x\in[0, \underline{s}_1(t)),~t> NT, \\[5pt]
\underline{u}_{x}(0, t)=0,~ \underline{u}(s_1(t),t)=0, ~ \underline{s}_1'(t)=-\mu_1\underline{u}_{x}(\underline{s}_1(t),t),~ t> NT,\\[5pt]
\underline{s}_1(NT)=s_1^{0}, ~\underline{u}(x, NT)=z(x, NT),~0\leq x<s_1^{0}.
\end{array}\right.
\end{align*}
In the case $0<s_{1}^{0}<s_{1}^{*}$, by Proposition 3.2 (ii), there exists $\overline{\mu}_1>0$ depending on $u_0$ and $s_1^0$ such that $\underline{s}_{1,\infty}=\infty$ provided that $\mu_1>\overline{\mu}_1$. Lemma 2.2 implies that $s_{1,\infty}\geq\underline{s}_{1,\infty}=\infty$ for $\mu_1>\overline{\mu}_1$. Therefore, when $0<s_{1}^{0}<s_{1,\min}\leq s_{1}^{*}$, we have
$s_{1,\infty}=\infty$ for $\mu_1>\overline{\mu}_1$.  \hfill $\Box$\\

\noindent\textbf{Lemma 4.10.}
$(i)$ Assume that $s_{1,\infty}=\infty$ and $s_{2,*}\leq s_2^0<\tilde{s}_{2,**}$, where
$\tilde{s}_{2,**}$ is the unique positive root of $\lambda_{1}(d_{2}, \alpha_{2}, b-hQ_{1},\cdot, T)=0$, then there exists $\hat{\mu}_2>0$ depending on $u_0, v_0,s_1^0, s_2^0$ and $\mu_1$ such that $s_{2,\infty}<\infty$ if $0<\mu_2\leq\hat{\mu}_2$.\\
$(ii)$ Assume that $s_{2,\infty}=\infty$ and $s_{1,*}\leq s_1^0<\tilde{s}_{1,**}$, where
$\tilde{s}_{1,**}$ is the unique positive root of $\lambda_{1}(d_{1}, \alpha_{1}, a-kQ_{2},\cdot, T)=0$, then there exists $\hat{\mu}_1>0$ depending on $u_0, v_0,s_1^0, s_2^0$ and $\mu_2$ such that $s_{1,\infty}<\infty$ if $0<\mu_1\leq\hat{\mu}_1$.\\

\noindent\textbf{Proof.}
We only deal with $(i)$, since $(ii)$ can be proved in a similar way.

Since $s_{1,\infty}=\infty$, by Lemma 4.3 $(ii)$, we have $\liminf_{n\rightarrow\infty}u(x,t+nT)\geq Q_{1}(x,t)$ locally uniformly in
$[0, \infty)\times[0,T]$. Let $\tilde{s}_{2,**}^\varepsilon$ be the unique positive root of $\lambda_{1}(d_{2}, \alpha_{2}, b-h(Q_{1}-\varepsilon),\cdot, T)=0$,
then for any $l>\tilde{s}_{2,**}^\varepsilon$ and small $\varepsilon>0$, there exists an integer $N\gg 1$ such that $u(x, t+nT)\geq Q_{1}(x,t)-\varepsilon$ for all $(x,t)\in[0,l]\times[0,T]$ and $n\geq N$.
Since $s_2^0<\tilde{s}_{2,**}^\varepsilon$, by Lemma 4.1 ($0<s_{2}^{\prime}(t)\leq C_{4}\mu_{2}$) and Corollary 2.1, there exists $\tilde{\mu}_2>0$ small enough, such that $s_2^{\mu_1, \mu_2}(NT)<\tilde{s}_{2,**}^\varepsilon$ for all $\mu_2\in(0, \tilde{\mu}_2]$.

Next, we apply comparison principle to derive that $v^{\mu_1, \mu_2}$ vanishes eventually provided that $\mu_2$ is small enough.
For each $\mu_2\leq\tilde{\mu}_2$, due to (4.14), we have
\begin{align*}
v_{t}^{\mu_{1},\mu_{2}}\leq d_{2}v_{xx}^{\mu_{1},\mu_{2}}-\alpha_{2}v_{x}^{\mu_{1},\mu_{2}}
+v^{\mu_{1},\mu_{2}}(b(x,t)-h(x,t)(Q_{1}(x,t)-\varepsilon)-v^{\mu_{1},\mu_{2}})
\end{align*}
for all $x\in [0, \min\{s_{2}^{\mu_{1},\mu_{2}}(t), l\}]$ and $t\geq NT$.
Let $(\bar{v}, \bar{s}_{2})$ be the solution of
\begin{align*}
\left\{\begin{array}{l}
\bar{v}_t=d_{2}\bar{v}_{xx}-\alpha_{2}\bar{v}_x
+\bar{v}(b(x,t)-h(x,t)(Q_{1}(x,t)-\varepsilon)-\bar{v}),
\quad 0<x<\bar{s}_{2}(t),\quad t>NT,\\[5pt]
\bar{v}_{x}(0, t)=0, ~\bar{v}(\bar{s}_{2}(t),t)=0, ~\bar{s}_{2}'(t)=-\mu_{2} \bar{v}_{x}(\bar{s}_{2}(t), t),\quad t>NT,\\[5pt]
\bar{s}_{2}(NT)=\bar{s}_{2}^{NT}, ~\bar{v}(x, NT)=\bar{v}^{NT}(x),\quad 0\leq x\leq\bar{s}_{2}^{NT},
\end{array}\right.
\end{align*}
where initial data $(\bar{s}_{2}^{NT}, \bar{v}^{NT}(x))$ satisfies $\bar{s}_{2}^{NT}\in (s_2^{\mu_1, \tilde{\mu}_2}(NT), \tilde{s}_{2,**}^\varepsilon)$ and
$\bar{v}^{NT}(x)>v^{\mu_1, \tilde{\mu}_2}(x, NT)$ for $0< x<s_2^{\mu_1, \tilde{\mu}_2}(NT)$.
Since $\bar{s}_{2}^{NT}<\tilde{s}_{2,**}^\varepsilon$, by Proposition 3.2 $(ii)$, there exists $\mu_{2}^{*}$ such that $\bar{s}_{2,\infty}\leq \tilde{s}_{2,**}^\varepsilon$ and
$\lim_{t\rightarrow \infty}\|\bar{v}(\cdot, t)\|_{C([0, \bar{s}_{2}(t)]}=0$
if $\mu_{2}\leq \mu_{2}^{*}$. By Corollary 2.1, we have
$s_{2}^{\mu_{1},\mu_{2}}(NT)\leq s_{2}^{\mu_{1}, \tilde{\mu}_{2}}(NT)<\bar{s}_{2}^{NT}$
and $v^{\mu_{1},\mu_{2}}(x, NT)\leq v^{\mu_{1},\tilde{\mu}_{2}}(x, NT)<\bar{v}^{NT}(x)$
for any $\mu_{2}\in(0, \min\{\tilde{\mu}_{2}, \mu_{2}^{*}\}]$
and $x\in [0, s_{2}^{\mu_{1},\mu_{2}}(NT)]$.
By Lemma 2.2, we can derive that $s_{2, \infty}^{\mu_1, \mu_2}<\bar{s}_{2,\infty}\leq \tilde{s}_{2,**}^\varepsilon<\infty$, which completes the proof of $(i)$. \hfill $\Box$\\

Let us discuss the effect of the coefficients $\mu_i$ $(i=1,2)$ on the spreading and vanishing.\\

\noindent\textbf{Theorem 4.4.}
Let $d_i$ $(i=1,2)$, $a, b, k$ and $h$ be fixed. Then\\
$(i)$ if $s_1^0<s_{1,min}$, then there exists $\mu_1^{**}\in[0, \bar{\mu}_1]$ depending on $u_0, v_0, s_1^0, s_2^0$ and $\mu_2$ such that $u$ spreads successfully if $\mu_1>\mu_1^{**}$ while $u$ vanishes eventually if $0<\mu_1\leq\mu_1^{**}$. Moreover, the following hold:

$(a)$ $\mu_1^{**}>0$ if one of the following is satisfied

\quad $(a.1)$ $0<s_1^0<s_{1,*}$;

\quad $(a.2)$ $s_{1,*}\leq s_1^0<\tilde{s}_{1,**}, 0<s_2^0<s_{2,min}$ and $\mu_2>\bar{\mu}_2$;

\quad $(a.3)$ $s_{1,*}\leq s_1^0<\tilde{s}_{1,**}$ and $s_2^0\geq s_{2, min}$.

$(b)$ $\mu_{1}^{**}=0$ if $s_1^0\geq s_{1,*}$, $0<s_2^0<s_{2,*}$ and $\mu_2\leq\underline{\mu}_2$.\\
$(ii)$ if $s_2^0<s_{2,min}$, then there exists $\mu_2^{**}\in[0, \bar{\mu}_2]$ depending on $u_0, v_0, s_1^0, s_2^0$ and $\mu_1$ such that $v$ spreads successfully if $\mu_2>\mu_2^{**}$ while $v$ vanishes eventually if $0<\mu_2\leq\mu_2^{**}$. Moreover, the following hold:

$(a)$ $\mu_2^{**}>0$ if one of the following is satisfied

\quad $(a.1)$ $0<s_2^0<s_{2,*}$;

\quad $(a.2)$ $s_{2,*}\leq s_2^0<\tilde{s}_{2,**}$, $0<s_1^0<s_{1,min}$ and $\mu_1>\bar{\mu}_1$;

\quad $(a.3)$ $s_{2,*}\leq s_2^0<\tilde{s}_{2,**}$ and $s_1^0\geq s_{1, min}$.

$(b)$ $\mu_{2}^{**}=0$ if $s_2^0\geq s_{2,*}$, $0<s_1^0<s_{1,*}$ and $\mu_1\leq\underline{\mu}_1$.\\

\noindent\textbf{Proof.}
We only deal with $(i)$, since $(ii)$ can be proved in a similar way.

Let $s_{1}^{0}<s_{1, min}$ be fixed. We defined
$A:=\{\mu>0~|~ s_{1,\infty}^{\mu, \mu_{2}}>s_{1, min}\}$. By Lemma 4.9 $(ii)$,
we know that $\mu_{1}^{**}:= \inf A\in [0, \bar{\mu}_{1}]$ is well defined.
Similar to the proof of Theorem 3 in \cite{chw15}, we can show that
$\mu_{1}^{**}\notin A$ if $\mu_{1}^{**}>0$.

We fist prove the part $(a)$. If $(a.1)$ holds, then Lemma 4.9 $(i)$ implies $\mu_{1}^{**}\geq \underline{\mu}_{1}>0$. If $(a.2)$ holds, by Lemma 4.9 $(iv)$ and Lemma 4.10 $(ii)$,
we have $\mu_{1}^{**}\geq \hat{\mu}_{1}>0$. If $(a.3)$ holds, then $s_{2,\infty}=\infty$, and
Lemma 4.10 $(ii)$ implies $\mu_{1}^{**}\geq \hat{\mu}_{1}>0$.

Next, we show the part $(b)$. Since $0<s_2^0<s_{2,*}$ and $\mu_2\leq\underline{\mu}_2$,
by Lemma 4.9 $(iii)$ and Theorem 4.1, we have
$\|v(\cdot, t)\|_{C[0, s_{2}(t)]}\rightarrow 0$ as $t\rightarrow \infty$. Similar to the proof
of Lemma 4.6, we can apply the assumption $s_1^0\geq s_{1,*}$ to prove that $s_{1,\infty}=\infty$, regardless of $\mu_{1}$. Thus, $\mu_{1}^{**}=0$, which completes the proof of $(i)$. \hfill $\Box$

\label{}








\end{document}